\documentclass[final,authoryear,5p,times,twocolumn]{elsarticle} 

\usepackage{lineno}
\modulolinenumbers[5]

\usepackage{graphicx}      
\usepackage{enumerate}

\usepackage{subfig}

\usepackage{url}
\usepackage{epstopdf}

\usepackage{amsmath}
\usepackage{amssymb}
\newtheorem{mydef}{Definition}
\newtheorem{myProp}{Proposition}

\journal{}

\biboptions{authoryear}

\newcommand{\set}[1]{\left\{#1\right\}}

\newcommand{\Real}{\mathbb R}

\begin{document}	

\begin{frontmatter}


\title{Risk minimization in life-cycle oil production optimization
} 




\author{Andrea Capolei}
\ead{acap@dtu.dk}

\author{Lasse Hjuler Christiansen} 
\ead{lhch@dtu.dk}

\author{John Bagterp J{\o}rgensen\corref{JBJO}}
\ead{jbjo@dtu.dk}

\address{Department of Applied Mathematics and Computer Science \& Center for Energy Resources Engineering, \\ Technical University of Denmark, DK-2800 Kgs. Lyngby, Denmark}

\cortext[JBJO]{Corresponding author.}






\begin{abstract}                
The geology of oil reservoirs is largely unknown. Consequently, the reservoir models used for production optimization are subject to significant uncertainty. To minimize the associated risk, the oil literature has mainly used ensemble-based methods to optimize sample estimated risk measures of net present value (NPV). However, the success in reducing risk critically depends on the choice of risk measure. As a systematic approach to risk mitigation in production optimization, this paper characterizes proper risk measures by the axioms of coherence and aversion. As an example of a proper measure, we consider conditional value-at-risk, $\text{CVaR}_{\alpha}$, at different risk levels, $\alpha$.  The potential of $\text{CVaR}_{\alpha}$ to minimize profit loss is demonstrated by a simulated case study. The case study compares $\text{CVaR}_{\alpha}$ to real-world best practices, represented by reactive control. It shows that for any risk level, $\alpha$, we may find an optimized strategy that provides lower risk than reactive control. However, despite overall lower risk, we see that all optimized strategies still yield some unacceptable low profit realizations relative to reactive control. To remedy this, we introduce a risk mitigation method based on the \textsl{NPV offset distribution}. Unlike existing methods of the oil literature, the offset risk mitigation approach minimizes the risk \textsl{relative} to a reference strategy representing common real-life practices, e.g. reactive control. In the simulated case study, we minimize the worst case profit offset to reduce the risk of realizations that do worse than the reactive strategy. The results suggest that it may be more relevant to consider the NPV offset distribution than the NPV distribution when minimizing risk in production optimization.
\end{abstract}


\end{frontmatter}



\section{Introduction}
Life-cycle production optimization seeks to enhance the process of oil recovery by maximizing a financial measure such as the net present value (NPV) over the expected reservoir life.  Simulation studies have demonstrated a significant potential of production optimization to increase overall profit. However, real-life applications are challenged by a wide range of uncertainties tied to reservoir simulation. To address the challenges of uncertainty, the oil literature has mainly considered \textsl{ensemble-based} methods. Such methods represent the uncertainty by a finite number of possible outcomes, i.e., by an \textsl{ensemble of realizations}. To minimize risk, the ensemble members are combined to form a sample estimated risk measure, which is optimized over the expected reservoir life. 

 The first ensemble-based approach, introduced by \cite{vanEssen:2009}, is known as robust optimization (RO). RO aims to maximize the life-cycle sample estimated \textit{expected return}. However, as discussed in \cite{Capolei:etal:2015b}, the expected profit is a risk neutral measure and it thereby neglects important risk indicators such as the lowest profit outcome. To account for risk directly, it has been proposed to include an additional measure, the profit standard deviation, to minimize the expected profit and reduce the profit standard deviation simultaneously. This approach is referred to as mean-variance optimization (MVO). Using this approach,  \cite{Yeten:2003}, \cite{Bailey:2005}, \cite{Alhuthali:2008}, and \cite{Capolei:etal:2014} were able to optimally trade-off expected profit and risk. However, as demonstrated by \cite{Capolei:etal:2015b}, the profit standard deviation may be misleading as a risk measure, when the profit distribution is asymmetrical. To remedy this, they suggest to measure the risk by coherent and averse risk measures such as conditional value-at-risk (CVaR). \textit{Conceptually}, $\text{CVaR}_\alpha$ represents the mean of the $\alpha$-percent lowest profit outcomes. The optimal operating strategy is the one that \textit{maximizes} the mean of the $\alpha$-percent lowest profit outcomes, i.e. provides the highest mean of the $\alpha$-percent lowest profit outcomes (see Fig. \ref{fig:procedureIdea}). \textit{Mathematically}, $\text{CVaR}_{\alpha}$ is defined as a risk measure, and the associated optimal strategy is defined as the strategy that \textit{minimizes} $\text{CVaR}_{\alpha}$ \citep{Rockafellar:2007}. Accordingly, $\text{CVaR}_{\alpha}$ of a profit distribution is connected to the negative of the profit distribution and is obtained as the mean of the $\alpha$-percent highest outcomes of the negative profit distribution. This value, $\text{CVaR}_{\alpha}$, is minimized. In the oil literature, CVaR has been used by \cite{Valladao:2013} as a deviation measure. They apply a  weighted sum method to find optimal trade-offs between expected profit  and profit deviation. They use $\text{E}_{\theta}(\psi)$ as a measure of profit and $\mathcal{D} = \text{E}_{\theta}(\psi) + \text{CVaR}_{\alpha}(\psi)$ as a deviation measure. As pointed out by  \cite{Capolei:etal:2015b}, the deviation measure, $\mathcal{D}$,  does not satisfy the axioms of risk aversion and monotonicity. This questions the ability of the measure to quantify risk. Recently, \cite{Siraj:etal:2015} and \cite{codas:etal:2015} used CVaR to minimize the risk of profit loss.
 
As a fundamental criterion to characterize proper risk measures for production optimization, this paper proposes to follow the axioms of coherence and aversion \citep{Rockafellar:2007}. As an example of a proper measure, we apply $\text{CVaR}_{\alpha}$ for different risk levels, $\alpha,$ to minimize the risk of profit loss. The potential of the method is illustrated by a numerical case study. The case study first minimizes the $\text{CVaR}_{\alpha}$  of the NPV distribution. All optimized control strategies are compared to a reference strategy that mimics a conventional reactive strategy in which producers are closed based on water breakthrough. We find that all optimized strategies provide lower measures of risk relative to the reference strategy. Secondly, the case study minimizes $\text{CVaR}_{\alpha}$ of the NPV offset distribution. The NPV offset distribution is computed as the difference between the fixed NPV distribution generated by the reference control strategy and the NPV distribution of the optimized strategy. 
 The profit offset distribution allows us to measure the probability that an optimized control strategy yields a lower profit than the reference strategy. In addition, we may measure the expected value of the profit difference. The results, for this case study, show that only for large values of the risk level, $\alpha$, the optimized control strategies have a lower risk than the reference control strategy. In other words, despite an overall lower risk, all optimized strategies have a positive probability of getting a lower NPV relative to the reference control strategy. This poses a significant risk of unacceptable low profit realizations. Therefore, in the case study we present a method to compute a control strategy that aims at optimizing the offset worst-case value. In this way, we manage to significantly increase the offset worst-case scenario. This result suggests that it may be more relevant to consider the NPV offset distribution as compared to the NPV distribution when minimizing risk in production optimization. 
 
 The paper is organized as follows. Section \ref{problemFormulation} formulates the oil production optimization problem under uncertainty as a risk minimization problem. Section \ref{coherMeasRisk} describes the basic properties that we require from an appropriate risk measure. In  Section \ref{RiskMeasuresExamples}, the CVaR measure is related to conventional risk measures. Section \ref{smoothCVarApprox} introduces a smooth approximation of CVaR that is appropriate for gradient-based optimization. Numerical results are presented in Section \ref{casestudy} and conclusions are made in Section \ref{conclusions}. \ref{Sec:app} lists the nomenclature used in this paper.

\section{Optimization under uncertainty}\label{problemFormulation}
In oil production optimization, the profit, $\psi $, can be considered as a function of the control vector, $u \in \mathcal{U} \subset \mathbb{R}^{n_u}$ with $\mathcal{U}$ expressing linear decision constraints, and a parameter vector, $\theta$; i.e.
\begin{equation}
	\psi = \psi(u; \theta). 
\end{equation}
The vector $\theta \subset \mathbb{R}^m$ represents the reservoir permeability field, porosity and economical parameters, etc.  Under the assumption of known parameters, $\theta$, we can maximize $\psi$ by solving the following deterministic optimal control problem \citep{Brouwer:Jansen:2004,Sarma:etal:2005,navdal:2006,Foss:2011,Volcker:etal:2011,Capolei:etal:2013}:
\begin{equation}\label{DeterministicOptimalControlProblem}
	\max_{u \in \mathcal{U}} \,\, \psi(u, \theta).
\end{equation}
However, due to the noisy and sparse nature of seismic data, core samples, borehole logs, future oil prices and plant costs, the parameters, $\theta$, are often highly uncertain. 
Mathematically, we may account for the uncertainty by considering the parameters, $\theta$, as a random variable with an associated probability distribution and uncertainty space, $\Theta$. Due to the complexity of real oil reservoirs, we only have incomplete information about the uncertainty space, $\Theta$. Therefore, the traditional way of modeling the uncertainty in oil production problems is to consider a finite set of possible outcomes for the parameters \citep{Krokhmaletal:2011,vanEssen:2009,Capolei:etal:2013,Capolei:etal:2014}. In particular, $\Theta$ is approximated by a discrete space $\Theta_d := \{\theta_1, \theta_2, \ldots, \theta_{n_d}\}$. As a consequence, a control input, $u$, will give rise to a finite set of possible profit outcomes $\psi^1 = \psi(u,\theta_1), \ldots, \psi^{n_d} = \psi(u,\theta_{n_d})$, with probabilities $p_1, \ldots, p_{n_d}$, where $p_i = \text{Prob}[\theta = \theta_i] \in [0,1]$ and $\sum_{i=1}^{n_d} p_i = 1$. 

When the parameters are uncertain and the associated profit outcomes are uncertain, the deterministic optimization problem  (\ref{DeterministicOptimalControlProblem}) cannot be extended directly to the optimization problem
\begin{equation}\label{StochasticOptimalControlProblem}
	\max_{u \in \mathcal{U}}  \,\, \psi(u, \theta \in \Theta_d).
\end{equation}
The optimization problem (\ref{StochasticOptimalControlProblem}) is undefined, as $\psi(u,\theta \in \Theta_d) = \set{\psi(u,\theta_1), \ldots,  \psi(u,\theta_{n_d})}$ is a set and not a function. To obtain a well defined problem, we replace the random variable, $\psi$, with a functional $\mathcal{R}$ that maps the random variable, $\psi = \psi(u,\theta)$, into a scalar deterministic measure. In this way, we may reformulate the undefined optimization problem (\ref{StochasticOptimalControlProblem}) to a well-defined optimization problem,
\begin{equation}\label{StochasticOptimalControlProblemWellDef}
	\min_{u \in \mathcal{U} } \,\, \mathcal{R}(\psi(u, \theta \in \Theta_d )). 
\end{equation}
Note that we have switched to a minimization problem. This is in accordance with the idea of interpreting $\mathcal{R}$ as a risk measure that we want to minimize. Further, it is known that minimizing $\mathcal{R}$ is the same as maximizing $-\mathcal{R}$ \citep{nocedal:book:2000}. In more detail, $\mathcal{R}$ is a surrogate for the distribution of $\psi$, where different $\mathcal{R}$ expressions capture different aspects of the profit distribution. In the oil community, different measures, $\mathcal{R}$, have been proposed. The proposed measures for $\mathcal{R}$ include the expected profit, i.e. $\mathcal{R}(\psi) = - \text{E}_{\theta}(\psi)$ \citep{vanEssen:2009}, and the mean-variance measure, i.e. $\mathcal{R}(\psi) = - \big( \lambda \text{E}_{\theta}(\psi) - (1-\lambda) \sigma^2(\psi) \big)$ with $\lambda \in [0,1]$ \citep{Capolei:etal:2014}. In Section \ref{RiskMeasuresExamples}, we argue why none of these measures can be considered to be satisfactory. 

Finally, we stress that this paper only focuses on single objective optimization. We do not consider important aspects connected to multi-objective optimization, e.g. the trade-off between long term vs short term profit \citep{vanEssen:2011,Christiansen:Capolei:Jorgensen:2016}. However, we note that our analysis on risk measures can be extended to such cases. We can e.g. use the weighted sum method \citep{liu:2015} to trade-off a long term profit measure $\mathcal{R}(\psi^{long})$ versus a short term profit measure  $\mathcal{R}(\psi^{short})$  by solving 
\begin{equation}
	\min_{u \in \mathcal{U} } \; \lambda \; \mathcal{R}(\psi^{short}) + (1-\lambda) \; \mathcal{R}(\psi^{long})
\end{equation}
for different values of $\lambda \in [0,1] $.

\subsection{Reference control strategy}\label{problemFormulationRelat}
In complex control optimization one often focuses on improving a reference control strategy, $u_{ref},$ that represents real-world best practices. In such cases, as an alternative to minimizing the risk of the profit distribution, $\mathcal{R}(\psi(u, \theta))$, directly, it may be more relevant to minimize the risk of the profit offset distribution,
\begin{equation}\label{offProfDistr}
	\psi_{off}(u, \theta) = \psi(u, \theta) - \psi(u_{ref}, \theta). 
\end{equation}
Here $\psi_{off}$ represents the profit offset with respect to a given reference case. From the profit offset distribution, we can extract two distributions of interest: the tail profit offset distribution,
\begin{align}\label{lowerTailProfDistr}
\{\psi_{off}(u, \theta) | \psi_{off} < 0 \},
\end{align}
and the upper tail profit offset distribution,
\begin{align}\label{upperTailProfDistr}
\{\psi_{off}(u, \theta) | \psi_{off} \geq 0 \}.
\end{align}
These distributions represent, respectively, the distribution of the profit loss and the profit gain with respect to the reference profit. In Section \ref{casestudy} we investigate the profit offset distribution as a tool for risk mitigation. 

\subsection{Net present value computation}\label{NPVcomp}
NPV is used as a measure of profit, $\psi$. The discrete profit distribution is given by the profit outcomes $\psi^i = \psi(u, \theta_i), i = 1, \ldots, n_d$, where \citep{Capolei:etal:2013,Capolei:etal:2014}
\begin{equation}\label{prodOpt:explicCost}
\begin{split}
	 \psi^i &= \psi(u, \theta_i)  \\
		 &= \sum_{k=0}^{N-1} \frac{\Delta t_k}{(1+d)^{\frac{t_{k+1}}{\tau}}} \Bigg[
			\overbrace{ \sum_{j \in \mathcal{P}} r_{o}\, q_{o,j}\big(u_k,x_{k+1}(u,\theta_i)\big) }^{\text{value of produced oil}}
	\\ & \qquad \qquad \qquad \quad \, - \overbrace{ \sum_{j \in \mathcal{P}} r_{wP}\, q_{w,j}\big(u_k,x_{k+1}(u,\theta_i)\big) }^{\text{cost of separating produced water}}  
	\\ & \qquad \qquad \qquad \quad  \, - \overbrace{ \sum_{j \in \mathcal{I}} r_{wI}\, q_{j}\big(u_k, x_{k+1}(u,\theta_i)\big)}^{\text{cost of injecting water}}
	\Bigg].
\end{split}
\end{equation}
 Subscript, $k$, denotes quantities at time $t_k$ for $t_0 < t_1 < \ldots < t_N$. Superscript, $i$, denotes the scenario, i.e. the simulation with parameter vector $\theta_i$ for $i=1,2, \ldots, n_d$. Accordingly, the reservoir state vector is denoted by $x_k^i = x_k(u,\theta_i) = x(t_k; u, \theta_i)$. The set of decision variables, $u = \{u_k\}_{k=0}^{N-1}$, is a sequence of piecewise constant control vectors, $u_k$ for $k=0,1, \ldots, N-1$, such that $u(t) = u_k$ for $t_k \leq t < t_{k+1}$ and $k=0,1, \ldots, N-1$. The flows at producer wells are denoted by $q_{w,j} = q_{w,j}(u(t),x(t))$ and $q_{o,j} = q_{o,j}(u(t),x(t))$; they are the volumetric water and oil flow rates at producer well $j \in \mathcal{P}$. The volumetric water flow at injector $j \in \mathcal{I}$  is denoted $q_j = q_j(u(t),x(t))$.   $r_o$, $r_{wP}$, and $r_{wI}$ represent the oil price, the water separation cost and the water injection cost, respectively. $d$ is the discount factor, $\Delta t_k = t_{k+1}-t_k$ is the time interval, and $N$ is the number of control steps. 
 
  In this paper, the states, $x_k^i = x(t_k; u, \theta_i)$, are computed using the two-phase immiscible flow model \citep{aziz:1979,chen:2006,chen:2007,Volcker:etal:2009}. However, the methodology can also be applied to black-oil models and compositional reservoir models. The flow rates at the production wells, $q_{w,j}$ and $q_{o,j}$ for $j \in \mathcal{P}$, are computed using the Peaceman well model \citep{Peaceman::1983}. The flow rate at injection wells may also be computed by the Peaceman well model. In this paper, we assume that the injection wells are rate controlled, and the decision variables, $u(t)$, are the injection flow rates.

 The profit offset distribution, $\psi_{off},$  is given by the profit offset outcomes
\begin{equation}\label{prodOpt:profExplOff}
	 \psi_{off}^i = \psi_{off}(u, \theta_i) = \psi(u, \theta_i) - \psi(u_{ref}, \theta_i), \quad i = 1, \ldots, n_d,
\end{equation}
where  $\psi(u, \theta_i)$ and  $\psi(u_{ref}, \theta_i)$ are the NPV outcomes (\ref{prodOpt:explicCost}) computed by using control  $u$ and $u_{ref},$ respectively.

\subsection{Risk minimization procedure}
\begin{figure*}[!tb]
	\begin{center}
		\subfloat[Model based production optimization under uncertainty.]{
			\includegraphics[width=0.8\textwidth]{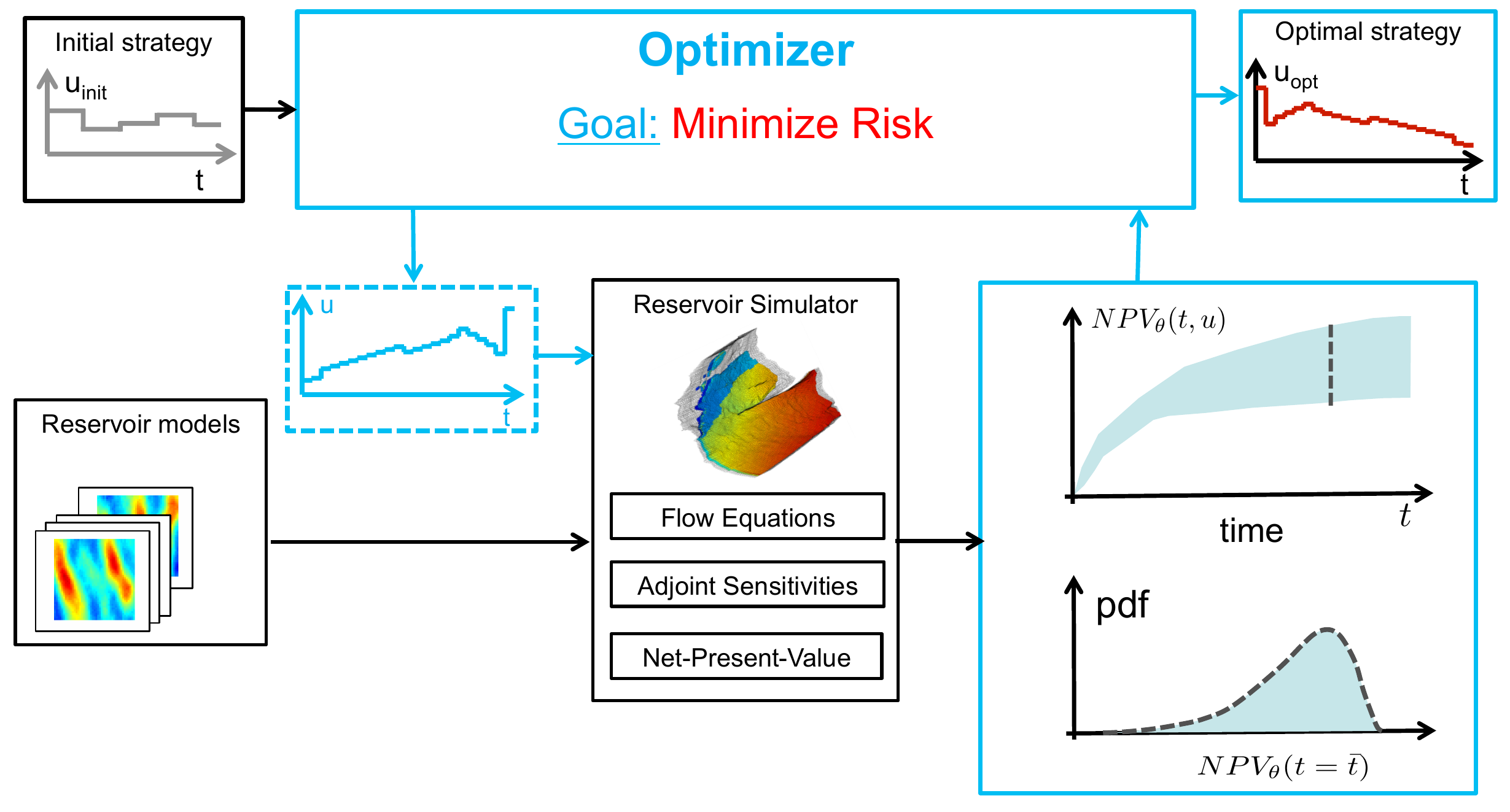}\label{fig:procedureIdea_Optim}			
		}	\\
		\subfloat[NPV risk mitigation. NPV uncertainty band versus time (left). NPV probability distribution function over the reservoir lifetime (right). Risk is reduced by maximizing the lifecycle average value (blue circles) of the $\alpha \%$ lowest profits (red areas). ]{
			\includegraphics[width=0.9\textwidth]{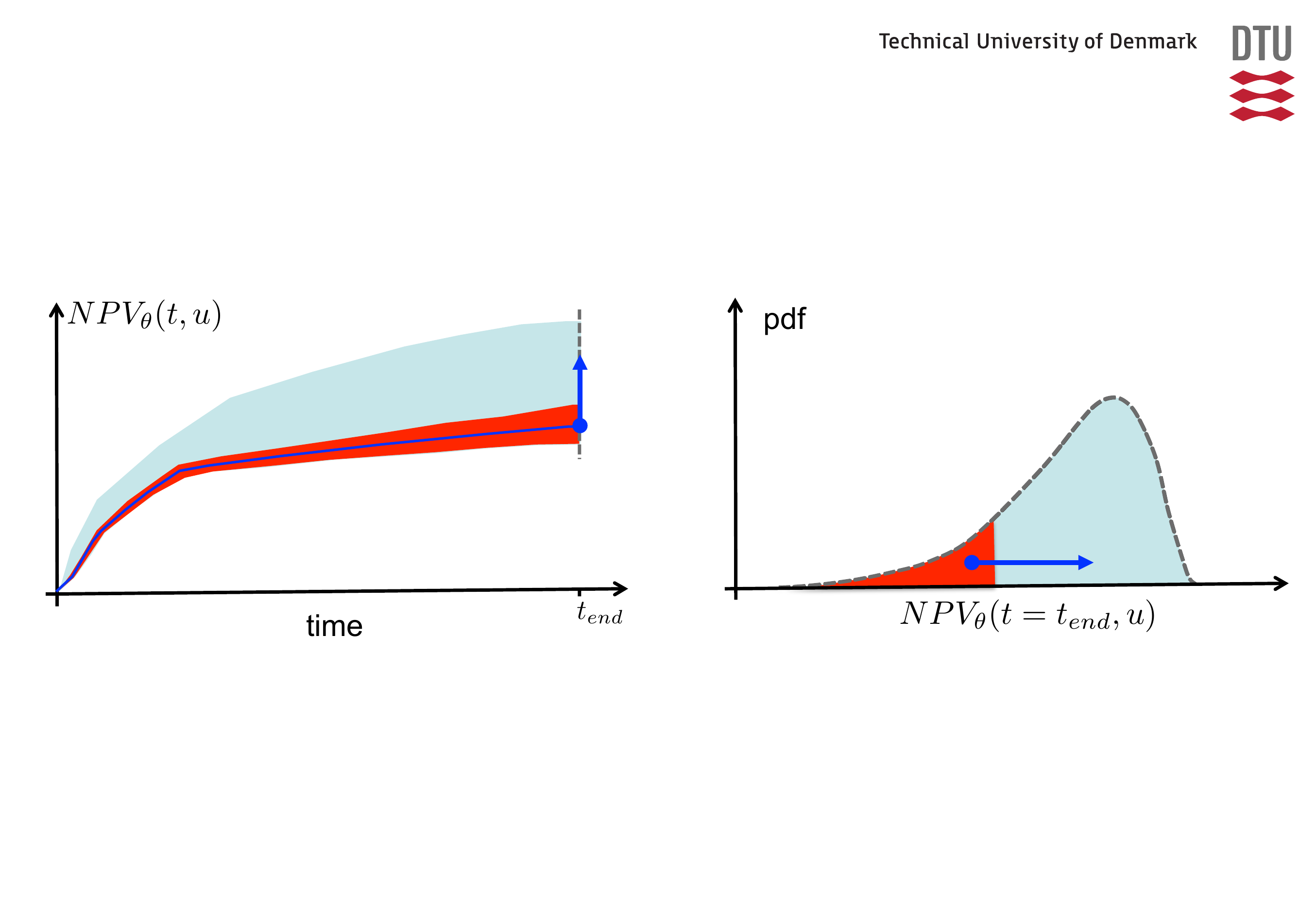}\label{fig:procedureIdea_NPV}			
		}	\\	
		\subfloat[NPV offset risk mitigation. NPV offset uncertainty band versus time (left). NPV offset probability distribution function over the reservoir lifetime (right).  Risk is reduced by maximizing the lifecycle average value (blue circles) of the $\alpha \%$ lowest offset profits (red areas).]{		
			\includegraphics[width=0.9\textwidth]{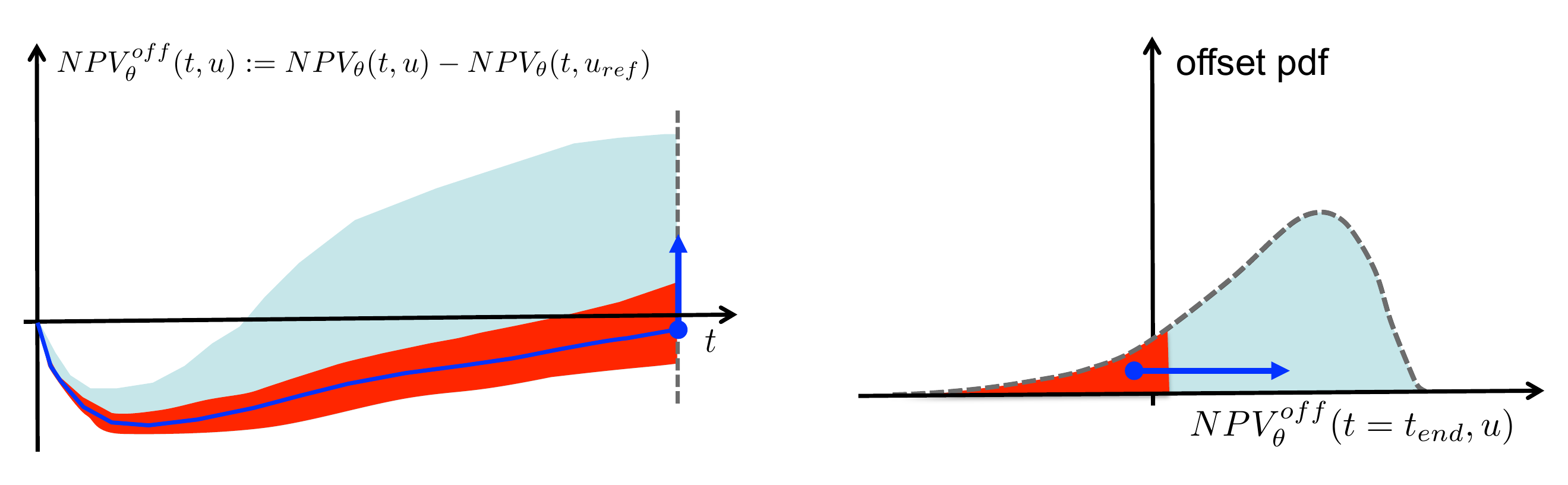}\label{fig:procedureIdea_NPVoff}
		}
	\end{center}
	\caption{Proposed risk minimization procedure. Model based optimization under uncertainty (\ref{fig:procedureIdea_Optim}) by the approaches of NPV risk minimization (\ref{fig:procedureIdea_NPV}) and NPV offset risk minimization (\ref{fig:procedureIdea_NPVoff}). The NPV risk mitigation approach minimizes risk by maximizing the lifecycle average value (blue circle) of the $\alpha \%$ lowest profits (red areas). The NPV offset strategy minimizes the risk relative to a reference strategy representing common real-life practices, e.g. reactive control. Unlike the NPV risk mitigation approach, the offset minimization method maximizes the lifecycle average value (blue circles) of the lowest $\alpha \%$  offset values of the profit distribution (red areas). The continuous blue curves in the figures represent the average of the $\alpha \%$ lowest profit values.  
	}\label{fig:procedureIdea}
\end{figure*}

Model based optimization seeks to determine the optimal control strategy that minimizes the risk of profit loss.
Fig. \ref{fig:procedureIdea_Optim} outlines the procedure, which consists of two key parts: 1) a reservoir simulator that, given an ensemble of reservoir models and a control input, computes the profit probability distribution. In Fig. \ref{fig:procedureIdea_Optim}, the black boxes and arrows mark this part of the procedure; 2) an optimizer that uses the profit distribution to compute an optimized control strategy that minimizes risk. This part of the procedure is marked by blue boxes and arrows in Fig. \ref{fig:procedureIdea_Optim}. The risk mitigation strategies proposed in this work, minimize risk by maximizing the lowest profit outcomes by two different approaches. The NPV risk mitigation approach maximizes the average value of the $\alpha$ percent lowest profits, whereas the approach of offset minimization maximizes the average value of the $\alpha$ percent lowest values of the offset  distribution. Fig. \ref{fig:procedureIdea_NPV} illustrates the concept of NPV optimization using the CVaR measure. Similarly, Fig. \ref{fig:procedureIdea_NPVoff} illustrates the use of CVaR for optimizing the NPV in relation to a reference strategy.

\section{Coherent averse measures of risk}\label{coherMeasRisk}
The role of a risk measure is to quantify the stochastic profit, $\psi$, by a numerical value, $\mathcal{R}(\psi),$ which serves as a surrogate for the overall profit distribution. This quantification of risk allows for fast and efficient decision-making. In particular, risk assessment of two scenarios, $\psi'$ and $\psi'',$ reduces to comparing the values $\mathcal{R}(\psi')$ and $\mathcal{R}(\psi'').$ However, the quality of the risk assessment critically depends on the properties of the risk measure in question. Therefore, it is important to have a characterization of properties that define a good risk measure.     
In this paper, we adhere to the coherence and aversion axioms introduced by \cite{Artzneretal:1999}, \cite{Rockafellar:2007}, and  \cite{Krokhmaletal:2011}.

Let $(\Omega, M, \text{Prob})$ be a probability space of elementary events, $\Omega$, with the sigma-algebra, $M$, over $\Omega$ and with a probability measure, Prob, on $(\Omega, M)$. Random profits are assumed to be measurable real-valued functions from $\mathcal{L}^2 (\Omega, M, \text{Prob})$ \citep{Zabarankinbook:2014}\footnote{$\mathcal{L}^2 (\Omega)$ is the Lebesgue space of measurable square-integrable functions on $\Omega$, i.e. $\psi \in \mathcal{L}^2(\Omega)$ is equivalent to $\int_{\Omega} | \psi(\omega) |^2 \; \text{dProb}[\omega] \; < \infty$.} and coherent averse risk measures are defined as:

\begin{mydef}\label{definitionCoAvRisk}
Coherent averse measures of risk are functionals $\mathcal{R}: \mathcal{L}^2 (\Omega) \rightarrow \mathbb{R}$ satisfying
\begin{enumerate}[{A}1.]
\item Risk aversion:
\begin{itemize}
 \item $\mathcal{R}(c) = -c$ for constants $c$ (constant equivalence)
 \item $\mathcal{R}(\psi)> -\text{E}_{\theta}(\psi)$ for  non-constant $\psi$ (aversion).
\end{itemize}
\item Positive homogeneity: $\mathcal{R}(\lambda \psi) = \lambda \mathcal{R}(\psi)$ when $\lambda>0$.
\item Sub-additivity: $\mathcal{R}(\psi'+\psi'') \leq \mathcal{R}(\psi') + \mathcal{R}(\psi'')$ for all $\psi'$ and $\psi''$.
\item Closure: $\forall c \in \mathbb{R}$, the set $\{ \psi | \mathcal{R}(\psi) \leq c  \}$ is closed. 
\item Monotonicity: $\mathcal{R}(\psi') \geq \mathcal{R}(\psi'')$ when $\psi' \leq \psi''$.
\end{enumerate}
\end{mydef}

These axioms require additional explanation. Axiom A1 formalizes the risk averse principle. A risk-averse decision maker does not rely on the expected profit exclusively and will always prefer a deterministic payoff of $\text{E}_{\theta}(\psi)$ over the stochastic payoff, $\psi$. 
The risk of a deterministic profit is given by its negative value, i.e. $\mathcal{R}(c) = -c$. This means that a specified deterministic profit value is less risky than a lower deterministic profit value. Further, it implies that $\mathcal{R}(\text{E}_{\theta}(\psi)) = -\text{E}_{\theta}(\psi)$, i.e. the condition $\mathcal{R}(\psi) > -\text{E}_{\theta}(\psi)$ can be restated as $\mathcal{R}(\psi) > \mathcal{R}(\text{E}_{\theta}(\psi))$ for $\psi \neq c$ and constant $c$. 
The positive homogeneity Axiom A2 ensures invariance under scaling, e.g. if the units of $\psi$ are converted from one currency to another, then the risk is simply scaled with the exchange rate. Finally, the positive homogeneity enables the units of measurements of $\mathcal{R}(\psi)$ to be the same as those of $\psi$. The sub-additivity Axiom A3 expresses the fundamental risk management principle of risk reduction via diversification. Further, combined with constant equivalence, i.e. $\mathcal{R}(c) = -c$, Axiom A3 ensures that
\begin{equation}\label{eq:translInvar}
	\mathcal{R}(\psi+c) = \mathcal{R}(\psi)-c.
\end{equation}
The property (\ref{eq:translInvar}) is referred to as \textsl{translational invariance}. This name relates to its financial interpretation. If $\psi$ is the payoff of a financial position, adding cash to this position reduces its risk by the same amount; in particular one has
\begin{equation}
	\mathcal{R}(\psi+ \mathcal{R}(\psi) ) = 0 = \mathcal{R}(0). 
\end{equation}
The translational invariance principle provides us with a natural way of quantifying acceptable risk \citep{Artzneretal:1999,Rockafellar:2007}. 
We consider the risk, $\mathcal{R}(\psi)$, acceptable if it is lower than the risk of obtaining a deterministic reference payoff $c_{ref}$, i.e. 
\begin{equation}\label{eq:lowerRefRisk}
	\mathcal{R}(\psi) \leq \mathcal{R}(c_{ref}) = -c_{ref}. 
\end{equation}
The monotonicity Axiom A5 implies that we consider $\psi'$ more risky than $\psi''$ if every realization of the profit $\psi'$ is smaller than every realization of the profit $\psi''$.  
In the literature, risk measures that satisfy axioms A1-A4 are called averse measures of risk \citep{Rockafellar:2007,Krokhmaletal:2011}. Risk measures that satisfy axioms A2-A5 and the constant equivalence property are called coherent risk measures in the sense of Artzner \citep{Artzneretal:1999,Krokhmaletal:2011}.
Finally, we note that positive homogeneity and the sub-additivity imply convexity of the risk measure $\mathcal{R}(\cdot)$ \citep{Rockafellar:2007,Krokhmaletal:2011}. The convexity property is important to risk minimization, since it allows the optimizer to find solutions that are globally optimal. Therefore, we would prefer convex production optimization problems. In oil optimization problems, however, $\psi = \psi(u,\theta)$  is non-convex with respect to the decision vector $u$. As a consequence,  $\mathcal{R}(\psi(u,\theta))$ is non-convex and the optimizer can only be expected to find a local minimizer. 

\subsection{Combining risk measures}\label{totCoherMeasRisk}
Risk measures which satisfy the coherence and aversion axioms can be combined to form new coherent and averse risk measures. In fact, we have the following propositions \citep{Rockafellar:2007}
\begin{myProp}\label{propCoAvRisk1}
If $\mathcal{R}_1, \ldots, \mathcal{R}_M$ are coherent averse measures of risk, then $\mathcal{R}$ defined as the convex combination of $\set{\mathcal{R}_m}_{m=1}^{M}$,
\begin{subequations}
\begin{alignat}{3}\label{eq:riskMeasComb1}
	& \mathcal{R}(\psi) = \sum_{m=1}^M \lambda_m \mathcal{R}_m (\psi), \\
	& \sum_{m=1}^M \lambda_m = 1, \\
	&  \lambda_m \geq 0,   \quad && m=1,2, \ldots, M,
\end{alignat}
\end{subequations}
is a coherent averse measure of risk.
\end{myProp}
\begin{myProp}\label{propCoAvRisk2}
If $\mathcal{R}_1$ is a coherent averse measure of risk, then
\begin{align}\label{eq:riskMeasComb2}
	\mathcal{R}(\psi) =  \lambda \mathcal{R}_1(\psi) - (1-\lambda) \text{E}_{\theta}(\psi), \quad \lambda \in (0,1], 
\end{align}
is a coherent averse measure of risk.
\end{myProp}
Propositions \ref{propCoAvRisk1} and \ref{propCoAvRisk2} allow us to define a coherent averse total risk measure $\mathcal{R}_{total}$ by
\begin{subequations}
 \begin{align}\label{eq:riskMeasComb3}
	& \mathcal{R}_{total}(\psi) =  \sum_{m=1}^{M} \lambda_m \mathcal{R}_m(\psi) - \lambda_{M+1}  \text{E}_{\theta}(\psi), \\
	& \sum_{m=1}^{M} \lambda_m > 0, \\
	& \lambda_{M+1} = 1 - \sum_{m=1}^{M} \lambda_m, \\
	& \lambda_m \geq 0, \quad m=1,2, \ldots, M, M+1,
\end{align}
\end{subequations}
where $\mathcal{R}_m$ for $m=1, \ldots, M$ are coherent and averse risk measures.

\section{Risk measures in production optimization}\label{RiskMeasuresExamples}
The following introduces the risk measures considered in this paper.  We present how each measure is computed numerically and comment on properties of coherence and aversion (Definition \ref{definitionCoAvRisk}). \cite{Capolei:etal:2015b} provide a more detailed review of traditional approaches to risk quantification in oil production optimization. We note that for ease of treatment, this paper assumes that the profit, $\psi$, is a discretely distributed random variable with an equiprobable probability distribution, $\text{Prob}[\psi = \psi^{i_k}] = p > 0$, $k=1, \ldots,n_d$, where $p = 1/n_d$ and ${i_k}$ is a set of indexes such that the profits are ordered, i.e. $\psi^{i_1} \leq \psi^{i_2} \leq \ldots \leq \psi^{i_{n_d}}$.
\subsection{Expected profit}\label{ex:Expected}
\begin{equation}\label{eq:ExpectedMes}
	\mathcal{R}\big(\psi(u,\theta)\big) = -\text{E}_{\theta}(\psi(u,\theta)).
\end{equation}
The expected return is a coherent measure of risk. It is widely used in oil production optimization, where it is referred to as RO \citep{vanEssen:2009,Capolei:etal:2013,Capolei:etal:2014}. Drawbacks of the method include risk neutrality. In particular, the method fails to address the paramount risks of extremely low profit realizations  \citep{Capolei:etal:2015b}. The expected profit is computed by
\begin{equation}
	\text{E}_{\theta}(\psi(u,\theta)) = \sum_{k=1}^{n_d} p \psi^{i_k}.
\end{equation}

\subsection{Worst-case scenario}\label{ex:Worst}
\begin{equation}\label{meas:worstCase}
	\mathcal{R}\big(\psi(u,\theta)\big) = -\inf_{\theta} \left( \psi(u,\theta) \right).
\end{equation}
The worst-case scenario measure is both coherent and averse. It has been investigated by e.g. \cite{Alhuthali:2010}. As a drawback, the measure does not take the probability distribution of $\psi$ into account. Consequently, the risk quantification is often too conservative. Using the ordered set of profits, i.e. $\set{\psi^{i_j}}_{j=1}^{n_d}$ with $\psi^{i_1} \leq \psi^{i_2} \leq \ldots \leq \psi^{i_{n_d}}$, the worst-case scenario is computed by
\begin{equation}\label{eq:WorsMes}
	-\inf_{\theta} \left( \psi(u,\theta) \right) = -\psi^{i_1}.
\end{equation}

\subsection{Value-at-risk ($\text{VaR}_{\alpha}$)}\label{ex:Var}
\begin{equation}
	\mathcal{R}\big(\psi(u,\theta)\big) = \text{VaR}_{\alpha}(\psi(u,\theta)).
\end{equation}
Value-at-risk, VaR, is one of the most widely used risk measures in
financial risk management, where it is a major competitor to the standard deviation measure \citep{jpmorgan:1994,jorion:1997}. Given a profit distribution, $\psi$, $\text{VaR}_{\alpha}(\psi)$ is defined as the negative $\alpha$-quantile
\begin{equation}
	\text{VaR}_{\alpha}(\psi) = - q_{\psi}(\alpha), \quad \alpha \in [0,1],
\end{equation}
where 
\begin{equation}
	q_{\psi}(\alpha) = \inf \, \{z \, \big| \, \text{Prob}[\psi \leq z] > \alpha \}.
\end{equation}
The quantile with $\alpha$ risk level, denoted by $q_{\psi}(\alpha)$, is defined as the lowest profit  value for which the probability that the profit, $\psi$, is lower or equal than $q_{\psi}(\alpha)$ is greater than $\alpha$. The VaR concept is closely related to probabilistic constraints (also called chance
constraints) that has been introduced by \cite{cooper:1958}. Probabilistic constraints have been widely used in disciplines such as operations research,
stochastic programming, systems reliability theory,
reliability-based design, and optimization based control \citep{ditlevsen:1996,Rockafellar:2010}. 
With a probabilstic constraint, we may declare that a profit, $\psi$, should exceed a certain predefined level, $c_{ref}$, with probability of at least $1-\alpha$, i.e. 
\begin{equation}\label{probConstraint}
	\text{Prob}[\psi \ge c_{ref}] \geq 1-\alpha, \qquad  \alpha \in [0, 1].
\end{equation}
In the case of $\alpha = 0$,  constraint (\ref{probConstraint})
reduces to the worst case approach described in Section \ref{ex:Worst}. In the case $\alpha=1$, the constraint (\ref{probConstraint}) is always satisfied and therefore without significance. From a risk reduction point of view, the probabilistic constraint (\ref{probConstraint}) has a dual aspect. One aspect is that for a fixed $\alpha$, we would like to find the highest value of $c_{ref}$ such that (\ref{probConstraint}) is satisfied. This ensures that with a probability greater than $1-\alpha$, the profit lower bound is the highest possible. On the other hand, for a fixed $c_{ref}$ value, we would like to have $\alpha$ as low as possible to increase the probability of having profits larger than $c_{ref}$. 
The probabilistic constraint (\ref{probConstraint}) is equivalent to
\begin{equation}\label{probConstraintEquiv}
	\text{Prob}[\psi < c_{ref}] \leq \alpha
\end{equation}
and it can be expressed as a constraint on the VaR of $\psi$ \citep{Krokhmaletal:2011,Zabarankinbook:2014}, i.e.
\begin{equation}\label{varConstraint}
	\text{VaR}_{\alpha} (\psi) \leq -c_{ref}.
\end{equation}
A major deficiency of $\text{VaR}_{\alpha}$ is that it does not take the tail of the profit distribution beyond the $\alpha-$ quantile level into account. As a consequence, it neglects risk of disastrous low profit realizations.
Further, $\text{VaR}_{\alpha}$ does not satisfy the sub-additivity Axiom A3 \citep{Artzneretal:1999}. 
Note also that $\text{VaR}_{\alpha}$ is discontinuous with respect to
the risk level, $\alpha$. This implies that a small change in the value
of $\alpha$ can lead to a significant jump in the risk estimate provided by $\text{VaR}_{\alpha}$.
The value-at-risk, $\text{VaR}_{\alpha}$, at risk level $\alpha$ of the discretely distributed stochastic profit, $\psi$, is given by
\begin{equation}
\label{eq:discrVar}
    \text{VaR}_{\alpha} (\psi) =
    \begin{cases}
        -\psi^{i_1}, &   \alpha \in [0, p), \\
        -\psi^{i_2}, &   \alpha \in [p, 2 p),\\        
		 \, \, \, \, \vdots & \, \, \vdots \\                        
		-\psi^{i_{n_d}}, &   \alpha \in [(n_d-1) p, 1]. \\
    \end{cases}
\end{equation}
This corresponds to the compact formulation
\begin{equation}
	\text{VaR}_{\alpha}(\psi) = 
	\begin{cases}
		-\psi^{i_{j+1}} & \alpha \in [jp, \, (j+1)p), \, j \in \set{0, 1,\ldots,n_d-1} \\
		-\psi^{i_{n_d}} & \alpha = 1
	\end{cases}
\end{equation}
Consequently, $\text{VaR}_{\alpha}(\psi)$ can be computed as
\begin{equation}
\text{VaR}_{\alpha}(\psi) = 
\begin{cases}
	-\psi^{i_{j+1}} & \alpha \in [0, \, 1), \, j = \text{floor}(\alpha/p),   \\
	-\psi^{i_{n_d}} & \alpha = 1
\end{cases}
\end{equation}
The function floor is a function that rounds down to the nearest integer.

\subsection{Conditional value-at-risk ($\text{CVaR}_{\alpha}$)}\label{ex:Cvar}
\begin{equation}
	\mathcal{R}(\psi(u,\theta)) = \text{CVaR}_{\alpha}(\psi(u,\theta)).
\end{equation}
 Conditional value-at-risk, CVaR, has been introduced to resolve the deficiencies of VaR \citep{Rockafellar:2002}. CVaR is defined as the average of VaR in the interval $[0,\alpha]$, i.e.
\begin{equation}
	\text{CVaR}_{\alpha}(\psi) = \frac{1}{\alpha} \int_0^{\alpha} \text{VaR}_s(\psi) ds, \quad \alpha \in [0, \,1].
\end{equation}
The CVaR measure is both coherent and averse when $\alpha \in (0,1)$.
As opposed to the VaR measure, CVaR is continuous with respect to the risk level, $\alpha$, it is sub-additive, and it accounts for the entire $\alpha$-tail of the distribution. 
The risk measure $\text{CVaR}_{\alpha}$ of a sorted discrete profit distribution, $\set{\psi^{i_j}}_{j=1}^{n_d} = \set{\psi(u,\theta_{i_j})}_{j=1}^{n_d}$, is given by
\begin{equation}
\label{eq:discrCVaR}
    \text{CVaR}_{\alpha}(\psi)   =
    \begin{cases}
        - \psi^{i_1},  & \alpha \in [0, p), \\
        - \frac{p \psi^{i_1} + (\alpha-p) \psi^{i_2}}{\alpha}, &  \alpha  \in [p, 2 p), \\  
        \quad \ldots &  \quad  \ldots  \\
         -\frac{ p ( \psi^{i_1}+\ldots+ \psi^{i_j}) + (\alpha - j p) \psi^{i_{j+1}} }{\alpha}, & \alpha \in [j p, (j+1) p), \\                    
                \quad \ldots  & \quad  \ldots  \\
		- \sum_{k=1}^{n_d} p \psi^{i_k},   &  \alpha  = 1. 
    \end{cases}
\end{equation}
Consequently, $\text{CVaR}_{\alpha}$ may be computed as
\begin{subequations}
	\label{eq:ComputationCVaRalpha}
		\begin{equation}
		\text{CVaR}_{\alpha}(\psi) = - \psi^{i_1}, \qquad \alpha \in [0, \, p),
	\end{equation}
	and
\begin{equation}
	\text{CVaR}_{\alpha}(\psi)
	=
	- \frac{1}{\alpha} 
	\left(
		\sum_{k=1}^{j} p \psi^{i_k} + (\alpha - jp) \psi^{i_{j+1}}
			\right), 
	\, \alpha \in [j p, \, (j+1)p),
\end{equation}
i.e. for $j=\set{1, \ldots, n_d-1}$ and $\alpha \in [p, \, 1)$. Given a value of $\alpha$, $j$ may be computed as $j = \text{floor}(\alpha/p)$. In the case $\alpha=1$, CVaR is computed by
\begin{equation}
	\text{CVaR}_{\alpha} (\psi)= - \sum_{k=1}^{n_d} p \psi^{i_k}, \quad \alpha=1.
\end{equation}
\end{subequations}
We note that $\text{CVaR}_{\alpha}$ reduces to the worst-case measure (\ref{eq:WorsMes}) for $\alpha \in [0,p)$, and $\text{CVaR}_{\alpha}$ reduces to the expected profit measure (\ref{eq:ExpectedMes}) for $\alpha = 1$. Further, if we consider the discrete risk levels $\alpha = j p$ for $j= 1, \ldots, n_d$,  (\ref{eq:ComputationCVaRalpha}) simplifies to 
\begin{equation}
    \text{CVaR}_{j p}(\psi)  = -\text{E}_{\theta}[\psi | \psi \leq \psi^{i_j}] = -\frac{\sum_{k=1}^{j} \psi^{i_k}}{j}, \, j = 1, \ldots, n_d.
\end{equation}
For these $\alpha$ values, $\text{CVaR}_{\alpha}(\psi) $ is simply the mean of the $j$ lowest profit outcomes.
Finally, we note that for any $\alpha \neq 1$, $\text{CVaR}_{\alpha}\big(\psi(u,\theta)\big)$ is non-differentiable with respect to the controls, $u$. Non-differentiable points occur when two or more of the profit outcomes have the same value \citep{Christiansen:Capolei:Jorgensen:2016}. In such cases the set of profit outcomes cannot be uniquely sorted in ascending order.

\section{Smooth approximation for $\text{CVaR}_{\alpha}$}\label{smoothCVarApprox}

The non-differentiability of the CVaR measure may interfere with the optimization procedure (\ref{StochasticOptimalControlProblemWellDef}) when $\mathcal{R} = \text{CVaR}_{\alpha}$ and $\alpha \in [0, \, 1)$. To improve convergence, we use an approximation of $\text{CVaR}_{\alpha}$ that is differentiable with respect to the control, $u$.
\cite{Rockafellaretall:2002} have shown that $\text{CVaR}_{\alpha}(\psi)$ for $\alpha \in (0,1)$ can be computed as the optimal value of the following optimization problem
\begin{equation}
\text{CVaR}_{\alpha}(\psi) = -  \max_{c \in \mathbb{R}} \left( c   - \frac{1}{\alpha} \text{E}_{\theta}\left( [c-\psi(u,\theta)]_+ \right) \right), \quad \alpha \in (0,1),
\end{equation}
where $[t]_+ = \max\{t,0\}$. As prooved by \cite{Zabarankinbook:2014}, the optimal $c$ is any value from the interval $[\text{VaR}_{1-\alpha}(-\psi), -\text{VaR}_{\alpha}(\psi)]$. Hence, as also shown by \cite{Zabarankinbook:2014},  the task of finding a control input, $u$, that minimizes $\text{CVaR}_{\alpha}$ is 
equivalent to solving the following optimization problem 
\begin{equation}\label{optProblemFull}
\min_{c \in \mathbb{R}, u \in \mathcal{U}} - c  + \frac{1}{\alpha} \text{E}_{\theta}\left( [c-\psi(u,\theta)]_+ \right), \quad \alpha \in (0,1).
\end{equation}
In a scenario approach, the expectation can be estimated by its sample average. Problem (\ref{optProblemFull}) is then approximated by
\begin{equation}\label{optProblemApproximate}
\min_{c \in \mathbb{R}, u \in \mathcal{U}}  -c  
+ \frac{p}{\alpha} \sum_{i=1}^{n_d} \max\{0, c-\psi(u,\theta_i) \},
\quad \alpha \in (0,1), 
\end{equation}
where we used the equiprobable assumption, $p_i = p = 1/n_d$ for the realizations. Even if the profit function $\psi(u,\theta_i)$ is continuously differentiable with respect to $u$ for all $\theta$, problem (\ref{optProblemApproximate}) is not directly tractable by common nonlinear optimization algorithms due to the nonsmoothness of the $\max$ function. To overcome this issue, the problem (\ref{optProblemApproximate})  is reformulated as the equivalent smooth problem
\begin{subequations}\label{optProblemApproximateSmooth}
	\begin{alignat}{5}
	 & \min_{c \in \mathbb{R}, u \in \mathcal{U}, y \in \mathbb{R}^{n_d} }  \; &&  
	 -c +\frac{p}{\alpha} \sum_{i=1}^{n_d} y_i, \\
	& s.t. &&  y_i \geq c - \psi(u, \theta_i), \quad  &&i=1,\ldots,n_d,  \\
	&      && y_i \geq 0,   &&i=1,\ldots,n_d,
	\end{alignat}
\end{subequations} 
where $\alpha \in (0,1)$.
 \cite{Rockafellar:2002} and \cite{Rockafellar:2010} showed that solving (\ref{optProblemApproximateSmooth}) is equivalent to minimize $\text{CVaR}_{\alpha}(\psi(u,\theta))$ on the control, $u$, i.e. $\min_{u \in \mathcal{U}} \, \text{CVaR}_{\alpha}(\psi(u,\theta))$. In this paper, we use the formulation (\ref{optProblemApproximateSmooth}) to minimize $\text{CVaR}_{\alpha}$, whenever $\alpha \in (0, \,1)$.	
\section{Case study}\label{casestudy}

In this section, we present a  numerical case study that uses a simulated reservoir to investigate $\text{CVaR}_{\alpha}$ as a tool for selecting an operating profile, $u_{opt}$, that minimizes the risk of obtaining low profits. 
In the first part of the case study, we minimize the $\text{CVaR}_{\alpha}$ of the NPV distribution (\ref{prodOpt:explicCost}) for 11 different risk levels, $\alpha$. The investigations include the optimization of the worst-case scenario, obtained by choosing $\alpha = 0.5\%$, and the optimization of the expected profit, obtained by choosing $\alpha = 100\%$. We use reactive control as the reference strategy.  The results show that for each optimized control strategy, $u_{opt}$, we obtain a lower risk than the risk associated with the reference strategy, $u_{ref}$, i.e.
\begin{equation*}
\text{CVaR}_{\alpha_l}\big(\psi(u_{opt}, \theta )\big) \leq \text{CVaR}_{\alpha_l}\big(\psi(u_{ref}, \theta )\big), \,\, l=1, \ldots, M=11.
\end{equation*}
%
In the second part of the case study, we minimize the risk of the NPV offset distribution, $\psi_{off}$. The NPV offset distribution, $\psi_{off}$, is computed by (\ref{prodOpt:profExplOff}).
 When using the NPV offset distribution, an optimized control strategy, $u_{opt},$ has a lower risk than a given reference strategy provided that
$
\text{CVaR}_{\alpha}\big(\psi_{off}(u_{opt}, \theta )\big) \leq 0.
$
It is important to note that $\text{CVaR}_{\alpha}(\psi_{off}(u_{opt},\theta)) \neq \text{CVaR}_{\alpha}(\psi(u_{opt},\theta)) - \text{CVaR}_{\alpha}(\psi(u_{ref},\theta)) $ where $\psi_{off} = \psi_{off}(u_{opt},\theta) =  \psi(u_{opt},\theta) - \psi(u_{ref},\theta)$.
The simulations in the case study show that when considering the NPV offset distribution, no strategy based on profit risks (as opposed to profit offset risks) has a lower risk  than the reference strategy for risk levels $\alpha < 20 \%$. Consequently, despite overall lower risk, all optimized control strategies have a positive probability of getting a lower NPV than the reative strategy. This is explained by significant negative offset worst-case profits, i.e. scenarios where the optimized strategies  perform worse than the reactive strategy. To improve this situation, we compute an optimized strategy that aims to maximize the offset worst-case value. In this way, we manage to significantly increase the worst case profit relative to the reactive strategy, which is the strategy used in current industrial practice and therefore the reference strategy.

\subsection{Reservoir model description}

\begin{figure}[!tb]
\begin{center}
\subfloat{
\includegraphics[width=0.5\textwidth]{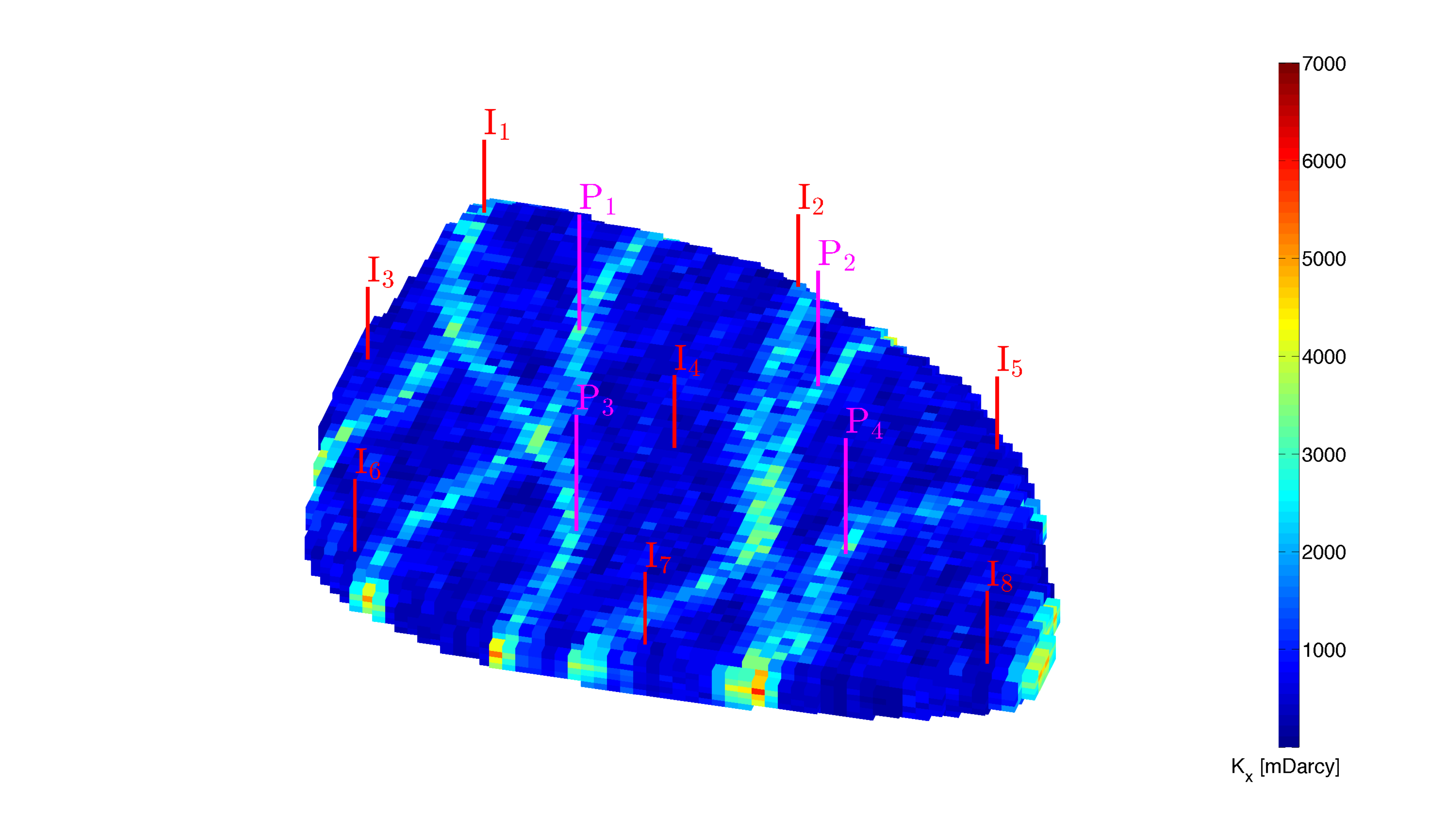}
}
\end{center}
\caption{Permeability field of Egg model for the first ensemble member used in the numerical simulations. The well configuration is also illustrated.}
 \label{fig:permFieldWells}
\end{figure}

The numerical simulations use the standard version of the Egg model \citep{eggmodel:2013}. This model has been used in a number of publications as a benchmark to test optimal control methodologies \citep{vanEssen:2009}. The Egg model is a synthetic reservoir model consisting of $60 \times 60 \times 7 = 25.200$ grid cells of which $18.553$ cells are active. The reservoir is produced for $3.600$ days under water flooding conditions. It contains eight water injectors and four producers, which are completed in all seven layers. The bhps of the producer wells are kept fixed at $395$ bar and the water injection rates are subject to control with a sample time of 90 days. The water injection rates are bound to be in the interval $[0, \, 79.5] \, \text{m}^3/\text{day}$. Fig. \ref{fig:permFieldWells} shows the well setup.
Model uncertainty is represented by an ensemble of 100 permeability realizations. 
Table \ref{tb:tableResDescr} provides petrophysical and economical simulation parameters. The flow in the reservoir is simulated using a two phase (oil and water) immiscible flow model with zero capillary pressure and incompressible fluids and rocks.

\begin{table}[!tb]
\begin{center}
\caption{Petro-physical and economical parameters for  the standard two phase Egg  model and the NPV objective function. 
} \label{tb:tableResDescr}
\scalebox{0.8}{
\begin{tabular}{llrr}
\hline\\[-8pt]
& Description                 & Value          & Unit \\\hline
$h$     & Grid-block height                    & $4$            & m \\
$\Delta x, \Delta y$     & Grid-block length/width                    & $8$            & m \\
$\phi$     & Porosity                    & $0.2$            & - \\
$c_o$      & Oil compressibility        & $1.0 \cdot 10^{-10}$              & Pa$^{-1}$ \\
$c_r$      & Rock compressibility        & $0$              & Pa$^{-1}$ \\
$c_w$      & Water compressibility        & $1.0 \cdot 10^{-10}$              & Pa$^{-1}$ \\
$\mu_o$    & Oil dynamic viscosity       & $5\cdot 10^{-3}$ & Pa $\cdot$ s \\
$\mu_w$    & Water dynamic viscosity     & $1.0\cdot 10^{-3}$ & Pa $\cdot$ s \\ \hline
$k_{ro}^{0}$     & End-point relative permeability, oil                     & $0.8$            & - \\
$k_{rw}^{0}$     & End-point relative permeability, water                     & $0.75$            & - \\
$n_o$      & Corey exponent, oil      & $4.0$            & -   \\
$n_w$      & Corey exponent, water    & $3.0$            & -   \\ \hline
$S_{or}$   & Residual oil saturation     & $0.1$            & -   \\
$S_{ow}$   & Connate water saturation    & $0.2$            & -   \\
$p_{c}$ & Capillary pressure  & $0$            & Pa   \\
$P_{init}$ & Initial reservoir pressure (top layer)  & $40 \cdot 10^{6}$            & Pa   \\
$S_{w,0}$     & Initial water saturation                    & $0.1$            & - \\
$p_{bhp}$ & Production well bottom hole pressures   & $39.5 \cdot 10^{6}$            & Pa   \\
$q_{wi,min}$ & Minimum water injection rate for well    & $0$            & m$^3$/day   \\
$q_{wi,max}$ & Maximum water injection rate for well    & $79.5$            & m$^3$/day   \\
$r_{well}$     & Well-bore radius                   & $0.1$            & m \\
$T$     & Simulation time                   & $3600$            & day \\
$N$     & Number of control steps                   & $40$            & - \\ \hline
$r_o$      & Oil price                   & $126$            & USD/m$^3$   \\ 
$r_{wP}$      & Water separation cost      & $19$             & USD/m$^3$   \\
$r_{wI}$      & Water injection cost       & $6$             & USD/m$^3$  \\ \hline
$d$      & Discount factor      &  			$0$              &    \\ \hline
\end{tabular} }
\end{center}
\end{table}

\subsection{Numerical optimization method}
To solve problem (\ref{StochasticOptimalControlProblemWellDef}), we use a gradient based optimization algorithm provided by Matlab's optimization toolbox  \citep{MATLAB}. Given an iterate of the optimizer, $\psi(u, \theta^i)$ is computed by solving the flow equations using MRST \citep{lieetall:2012}. The gradient, $\nabla_{u} \psi$ is computed by the adjoint method \citep{Jorgensen:2007,Volcker:etal:2011,Capolei:ECMOR:2012,Capolei:etal:2012,Jansen:2011,Sarma:etal:2005,suwartadi:2012}. An optimal solution is reported if the KKT conditions are satisfied to within a relative and absolute tolerance of $10^{-6}$. The current best but non-optimal iterate is returned in cases  for which the optimization algorithm uses more than 400 iterations, the relative change in the cost function is less than $10^{-6}$, or the relative change in the step size is less than $10^{-10}.$ These stopping criteria are independent, i.e. when one of the criteria is  satisfied, the optimizer stops. Furthermore, the cost function is normalized to improve convergence. The normalization consists of dividing by $10^6$ such that the objective function is appropriately scaled. For a given risk level, $\alpha$, we optimize $\text{CVaR}_{\alpha}$ by using three different initial guesses. Among the 3 optimized solutions, we select the optimized strategy that yields the lowest $\text{CVaR}_{\alpha}$ value. These initial guesses are constant water injection trajectories of $24$, $40$ and $60$ m$^3$/day.

\subsection{Profit risk minimization}\label{EggRiskMin}
The following demonstrates the potential of $\text{CVaR}_{\alpha}$ to minimize the risk of profit loss. In particular, we minimize the following $11$ risk measures 
\begin{equation}\label{riskMeasAcc}
	\text{CVaR}_{0.5 \%}(\psi(u,\theta)), \, \{\text{CVaR}_{(j \cdot 10) \%}(\psi(u,\theta))\}_{j=1}^9, \, -\text{E}_{\theta} \left( \psi(u,\theta) \right).
\end{equation}
Note that the RO strategy corresponds to the objective function $\text{CVaR}_{100\%}(\psi(u,\theta)) = - \text{E}_{\theta}\left( \psi(u,\theta)\right)$, and the worst case strategy corresponds to the objective function $\text{CVaR}_{0.5\%}(\psi(u,\theta)) \approx \lim_{\alpha \rightarrow 0} \text{CVaR}_{\alpha}(\psi(u,\theta)) = -\inf_{\theta} ( \psi(u,\theta))$. 
For each meausure used as the objective function, we solve problem (\ref{StochasticOptimalControlProblemWellDef}). Therefore, we obtain $11$ different optimized strategies, $u_{opt}$. The strategies are named according to the conventions provided in Table \ref{tb:namingConvContr}. As a representative of real-world best practices, we use the reactive control as reference strategy, $u_{ref}$. The reactive strategy is computed using a constant water injection rate of $60$ m$^3$/day. 

\begin{table}[!tb]
\begin{center}
\caption{Naming convention for the optimized control strategies.}\label{tb:namingConvContr}
\scalebox{0.8}{
\begin{tabular}{lll}
\hline\\[-8pt]
 Obj. function 
 & Control strategy name    & Abbreviation     \\ \hline
\vspace{-0.3cm}  \\ 
$\text{CVaR}_{0.5 \%}(\psi)$     			 & Worst case optimization            & w.c. opt          \\
$\text{CVaR}_{(j \cdot 10) \%}(\psi)$     	& Control strategy $(j \cdot 10) \%, j=1,\ldots, 9$            & c.s. $(j \cdot 10) \%$   \\
$-\text{E}_\theta (\psi)$     			 & Robust Optimization             & RO    \\ \vspace{-0.4cm}     \\ \hline 
\vspace{-0.3cm}  \\ 
$\text{CVaR}_{0.5 \%}(\psi_{off})$      			 & Offset worst case optimization              & off w.c. opt. \\ \vspace{-0.3cm}  \\ 
\hline
-     			 & Reference (reactive strategy)             & ref      
      \\ \hline
\end{tabular} }
\end{center}
\end{table}

Table \ref{tb:ExpSigWorst} compares the worst-case profit, the expected profit, and the profit standard deviation. We observe that all optimized strategies provide a higher expected profit than the reactive strategy. The improvements range from $2.4  \%$ to $3.2 \%$. The RO strategy has the highest expected profit. The expected profits of all the optimized strategies are in the same range. The largest relative difference is $0.9 \%$. Further, the profit standard deviations of the optimized strategies are comparable to that of the reference strategy. Notice also that all optimized strategies improve the worst-case profit compared to the reference strategy. As expected, the worst-case optimization strategy presents the highest worst-case profit value with an increased value of $5.8 \%$ relative to the reference strategy. Note also that the worst-case control strategy provides the lowest profit standard deviation. This comes at the price of a slightly lower expected profit compared to the RO strategy and the c.s. 30\% strategy.
\begin{table}[!tb]
	\begin{center}
		\caption{Key performance indicators for the NPV distributions of the profit risk minimization control strategies and the reactive strategy (ref).}\label{tb:ExpSigWorst}
		\scalebox{0.8}{
			\begin{tabular}{lccccccc}
				\hline\\[-8pt]
				Control 
				& $\inf (\psi)$    & $\text{E}_{\theta}(\psi)$ & $\sigma_{\theta}(\psi)$ & $5\%$ perc. & $95\%$ perc.  & -CVaR$_{30\%}$   \\
				strategy & 	$10^6$	 	 				& $10^6$	 & $10^6$ & $10^6$  & $10^6$  & $10^6$  \\
				& USD & USD & USD & USD & USD & USD \\ \hline
				w.c. opt.   		 &  42.94             &	45.43 &  1.43    & 43.34 & 47.83 & 43.70    \\
				c.s. $10 \%$   			 &  42.14        &  45.60 &  1.44 & 43.47 & 48.13 & 43.88         \\
				c.s. $20 \%$   			 &  41.94        &  45.62 &  1.44   & 43.47 & 48.12 & 43.89       \\
				c.s. $30 \%$   			 &  41.65        &  45.79 &  1.48  & 43.15 & 48.34 & 44.03        \\
				c.s. $40 \%$   			 &  41.70       &  45.76 &  1.49 & 43.13 & 48.30 & 43.98         \\
				c.s. $50 \%$   			 &  41.48       &  45.50 & 1.51 & 43.13 & 47.98 & 43.67         \\
				c.s. $60 \%$   			 &  41.49        &  45.56 &  1.53 & 43.10 & 48.11 & 43.72         \\
				c.s. $70 \%$   			 &  40.94       &  45.72 &  1.64 & 42.79 & 48.31 & 43.76         \\
				c.s. $80 \%$   			 &  41.46      &  45.75 &  1.55  & 43.01 & 48.26 & 43.90        \\
				c.s. $90 \%$   			 &  41.23        &  45.58 &  1.54 & 42.87 & 48.12 & 43.73         \\
				RO		   			 &  41.45    &  45.82 &  1.58 & 43.06 & 48.45 & 43.95         \\ \hline
				ref   			 &  40.60            &  44.38 &  1.57   & 41.59 & 46.57 & 42.46	       \\ \hline
			\end{tabular} 
		}
	\end{center}
\end{table}

\begin{figure*}[!tbh]
\begin{center}
\subfloat[Cumulative distribution function of the NPV.]{
\includegraphics[width=0.5\textwidth]{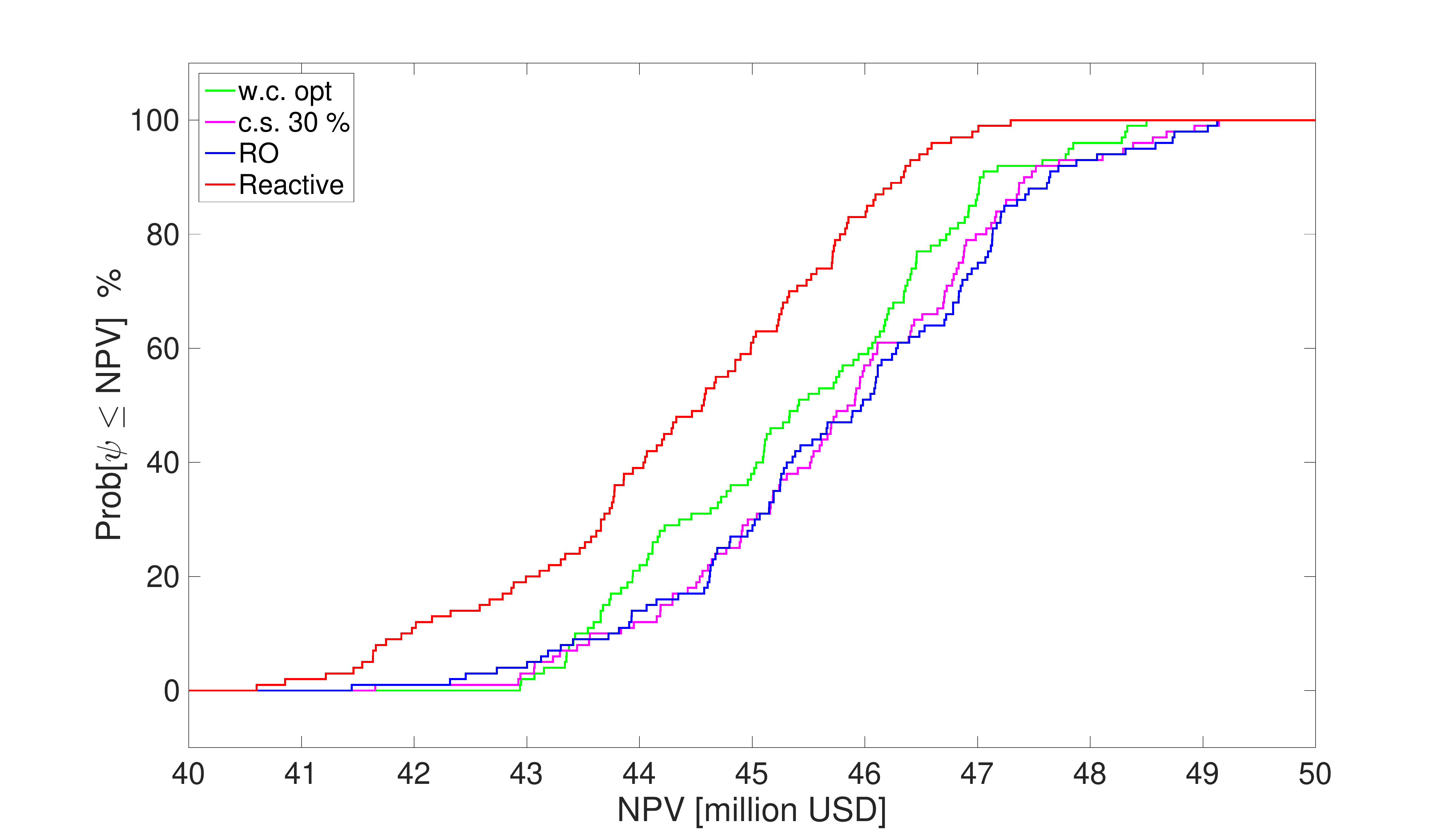}
\label{fig:cumDistr}} 
\subfloat[ $\text{CVaR}_{\alpha}(\psi)$ as a function of the risk level, $\alpha$.]{
\includegraphics[width=0.5\textwidth]{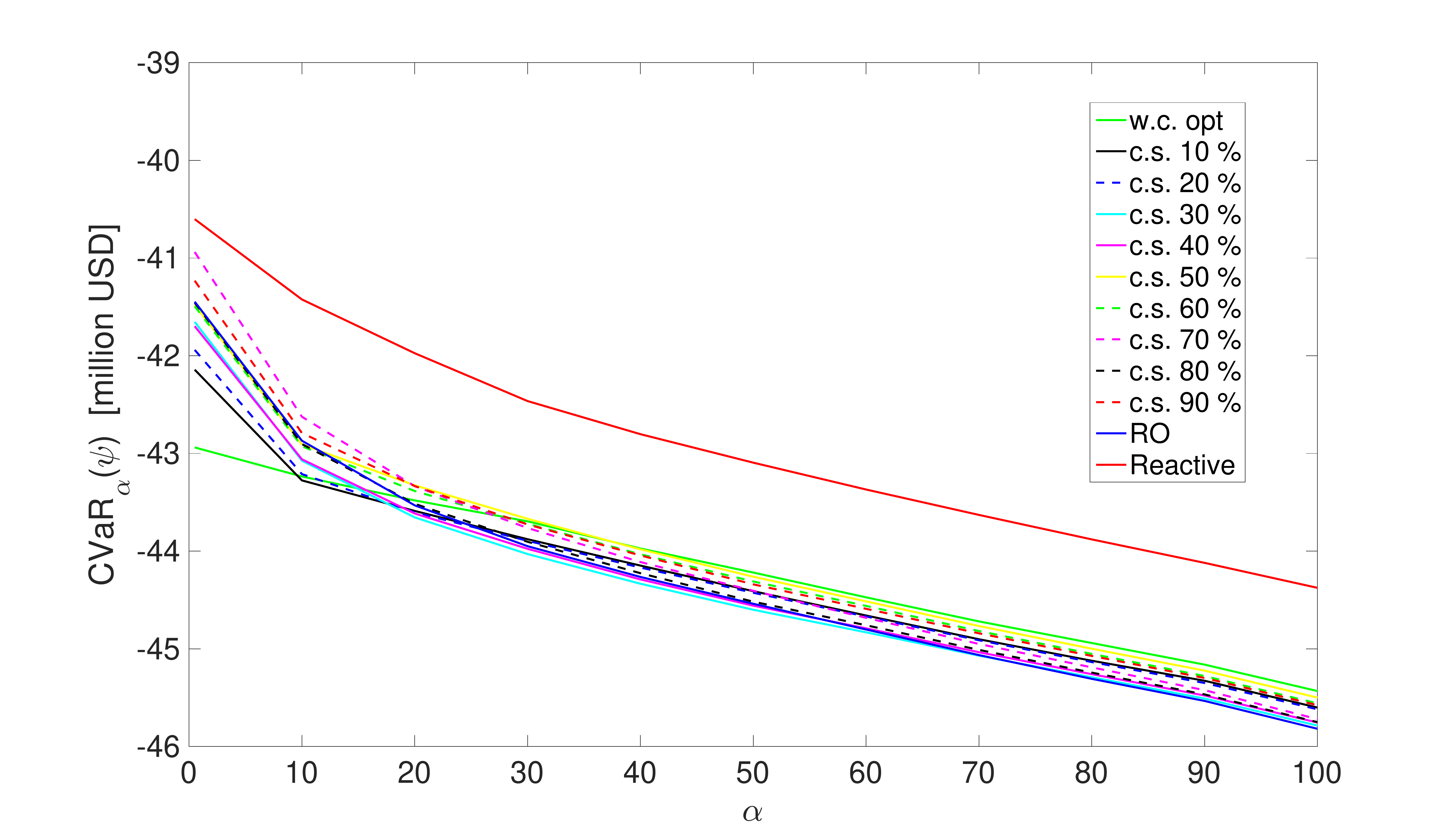}
\label{fig:cvarVsCvar}}
\end{center}
\caption{The probability function and the CVaR function of the NPV. (a) The cumulative distribution function for the NPV of selected control strategies. The reactive strategy has a higher risk than the optimizated stragies for a low profit realization. The RO control strategy has a higher risk of low NPV outcomes than the worst case control strategy and the c.s. 30\% control strategy. (b) The risk measure, CVaR$_{\alpha}$, for all control strategies as function of $\alpha$. For all values of $\alpha$, the risk of the reactive strategy is larger than the risk of the optimized strategies. The optimized strategies are similar for $\alpha > 20\%$. Large variations occur for $\alpha < 20\%$. This indicates that the main differences of the optimized strategies are related to how the low outcome realizations are distributed.}
 \label{fig:profPsiDistrAll}
\end{figure*}

Fig. \ref{fig:cumDistr} shows the cumulative distribution functions of selected strategies. Note that the reference strategy always has a higher probability of low profit realizations compared to the optimized strategies. This also holds for the optimized  strategies that are not included in the figure. These observations are supported by Fig. \ref{fig:cvarVsCvar}, which demonstrates that all optimized strategies provide lower risk than the reference control strategy. In accordance with (\ref{eq:lowerRefRisk}), all optimized strategies therefore have an acceptable risk. 

\begin{figure*}[!tbh]
\begin{center}
\subfloat[Strip chart of the NPV realizations for different control strategies.]{
\includegraphics[width=0.5\textwidth]{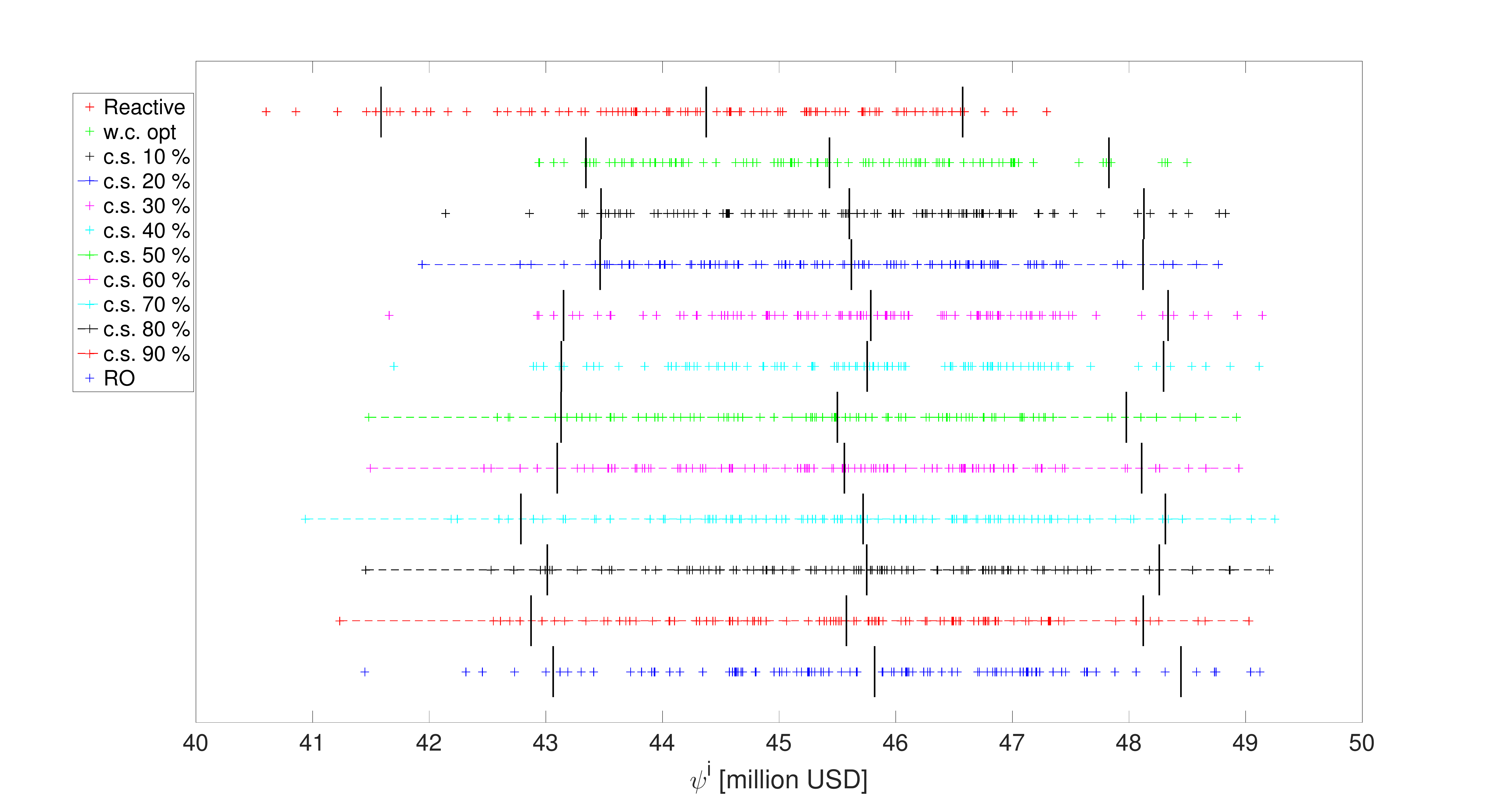}
 \label{fig:profPsiDistr}}
\subfloat[The total risk measure computed by (\ref{totRiskCase1}) for different control strategies, c.s. $\alpha\%$.]{
\includegraphics[width=0.5\textwidth]{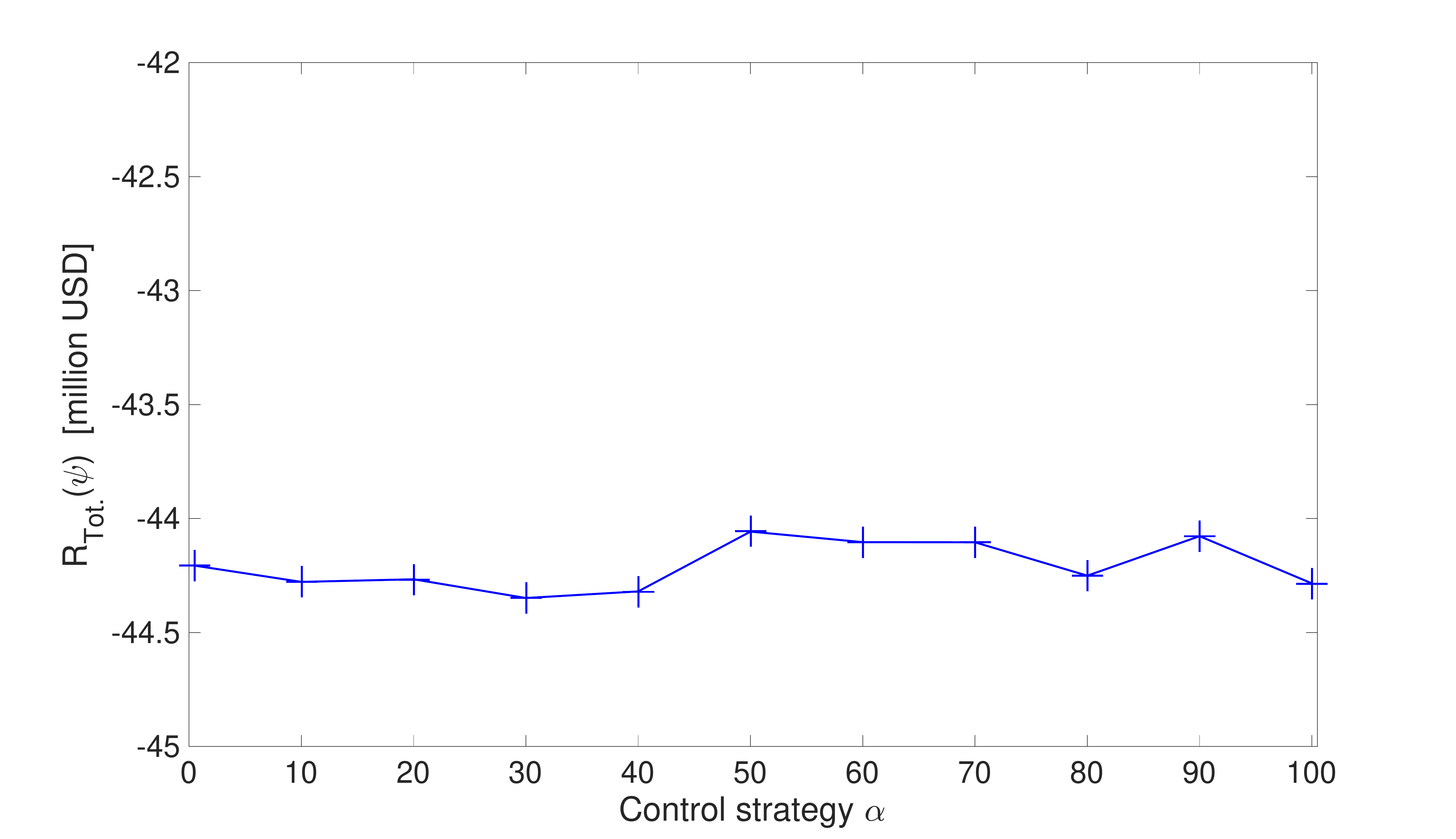}
\label{fig:totRiskCase1}}
\end{center}
\caption{The NPV realizations for selected control strategies and the overall risk measure. (a) Strip chart of the NPV realizations for different control strategies. The plot indicates the NPV of each control strategy for all 100 realizations. The black vertical lines indicate the 5\% perecentile, the mean, and the 95\% percentile. The reactive control strategy has the highest risk of low outcomes and also has the lowest mean NPV.  (b) The total risk measure computed by (\ref{totRiskCase1}) for $u_{opt}$ from different control strategies, i.e. c.s. $\alpha\%$. $\alpha = 0\%$ corresponds to the worst case control strategy and $\alpha = 100\%$ represents the RO control strategy. The total risk measure achieves its minimum for c.s. 30\%, i.e. $\alpha=30\%$. }
 \label{fig:contrSelec}
\end{figure*}

\begin{figure*}[!tbh]
\begin{center}
\includegraphics[width=0.995\textwidth]{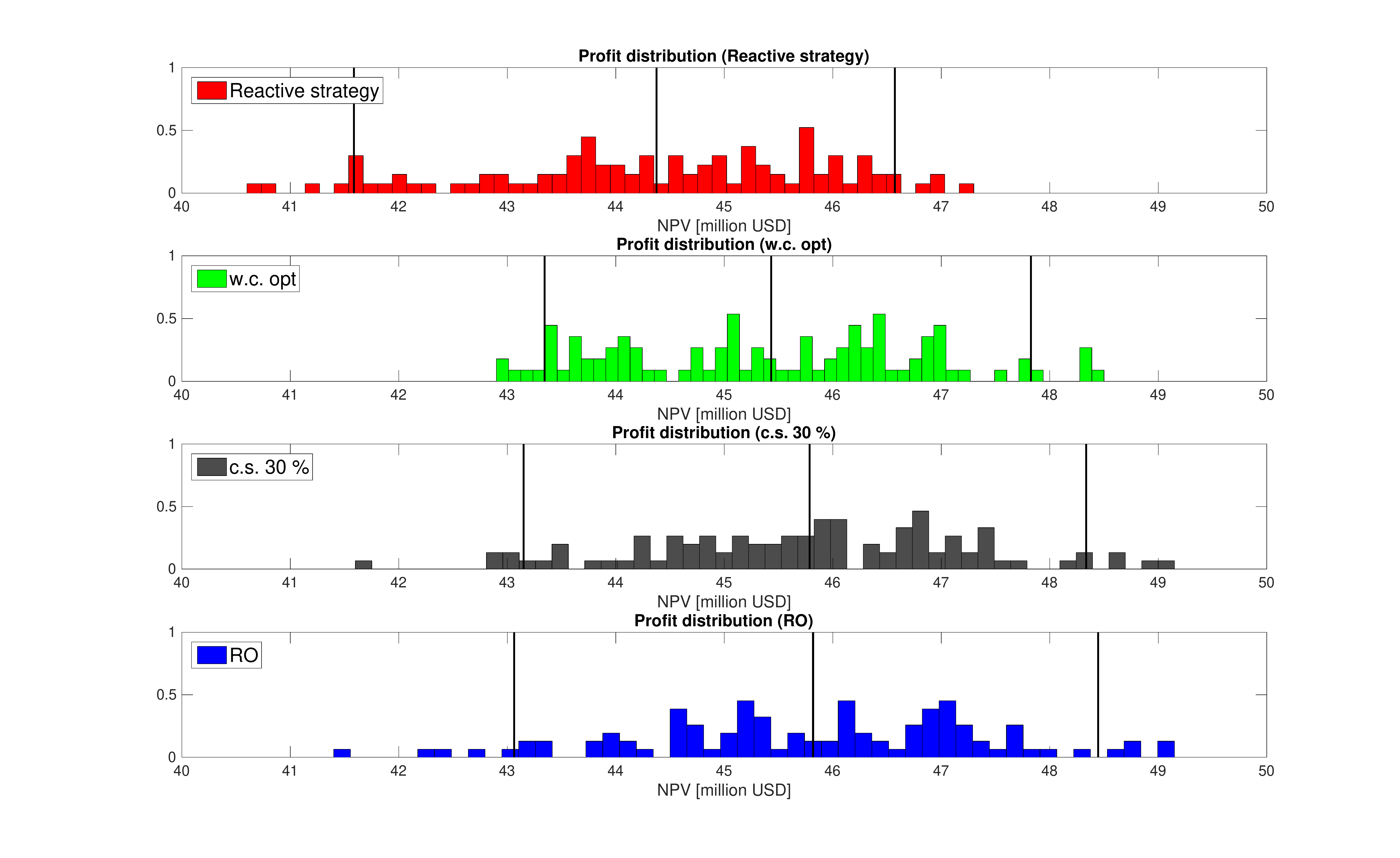}
\end{center}
\caption{Frequency plots of the profit distribution for selected control strategies. The black vertical lines indicate the 5\% percentile, the mean, and the 95\% percentile of the NPV distribution. The reactive control strategy has the lowest mean NPV and also the highest risk of low NPV outcomes. The RO control strategy has a mean NPV that is only slightly higher than the mean NPV of the c.s. 30\% control strategy. However, the RO control strategy has a higher risk of low NPV outcomes than the c.s. 30\% control strategy.}
 \label{fig:BarPlot}
\end{figure*}

Note from Fig. \ref{fig:cvarVsCvar} that while the values of $\text{CVaR}_{\alpha}$ are comparable for all optimized strategies whenever $\alpha >20 \%$,  large differences in $\text{CVaR}_{\alpha}$ occur for smaller choices of $\alpha.$  
This shows that the main differences between the optimized strategies are tied to the way the lowest profit realizations are distributed. This also follows from Fig. \ref{fig:profPsiDistr}, which presents strip charts of the profit distributions associated with the respective strategies. In particular, with the exception of the worst-case strategy and the reactive control, the optimized strategies have a significant gap between the worst-case profit and the remaining 99 profit realizations. As a way to quantify which optimized strategy has the overall lowest risk, we introduce the following total risk measure 
\begin{equation}\label{totRiskCase1}
	\mathcal{R}_{total}(\psi) = \frac{1}{11} \left(
	 \text{CVaR}_{0.5\%}(\psi) +  \sum_{j=1}^{10} \text{CVaR}_{j \cdot 0.1}(\psi) \right).  
\end{equation}
By using this total risk measure, we can identify the optimized strategy that provides the minimal average conditional value-at-risk over all the risk levels. Note that by Propositions \ref{propCoAvRisk1} and \ref{propCoAvRisk2}, the total risk measure is averse and coherent. Fig. \ref{fig:totRiskCase1} shows the total risk measure (\ref{totRiskCase1}) applied to all optimized strategies. We observe that c.s. 30\% has the lowest total risk. To illustrate the risk mitigation effects on the profit distribution, 
Fig. \ref{fig:BarPlot} compares c.s. 30\% to conventional strategies of RO and reactive control in term of a frequency plot. As indicated by the 5\% percentile, the risk of low profit realizations in the c.s. 30\% has been reduced noticeably compared to both the reactive control strategy and the RO control strategy.

%

\subsection{Analysis of the profit offset distribution and risk}
The results of the first part of the case study show that all optimized strategies provide lower risk related to low profit realizations than the reactive control. Therefore, from an overall perspective, the optimized strategies perform better than the reference strategy, i.e. than reactive control. Nevertheless, most of the optimized strategies give rise to worst case realizations with low profits. As is noticeable from Fig. \ref{fig:profPsiDistr} and Fig. \ref{fig:BarPlot}, the optimized strategies have lower overall risk of realizations with low NPV than the reactive strategy. However, these plots do not indicate if the low outcomes of the optimization based strategies and the reactive strategy occur for the same realizations. It may be that the optimized strategies despite overall lower risk contain a risk of yielding lower profits than the reactive strategy. Reservoir asset managers may be more concerned about the risk of doing worse than current best practice, i.e. using the reactive control strategy, than in the overall risk.

\begin{figure*}[!tb]
	\begin{center}
		\subfloat[Profit offset realizations using the worst case control strategy.]{
			\includegraphics[width=0.495\textwidth]{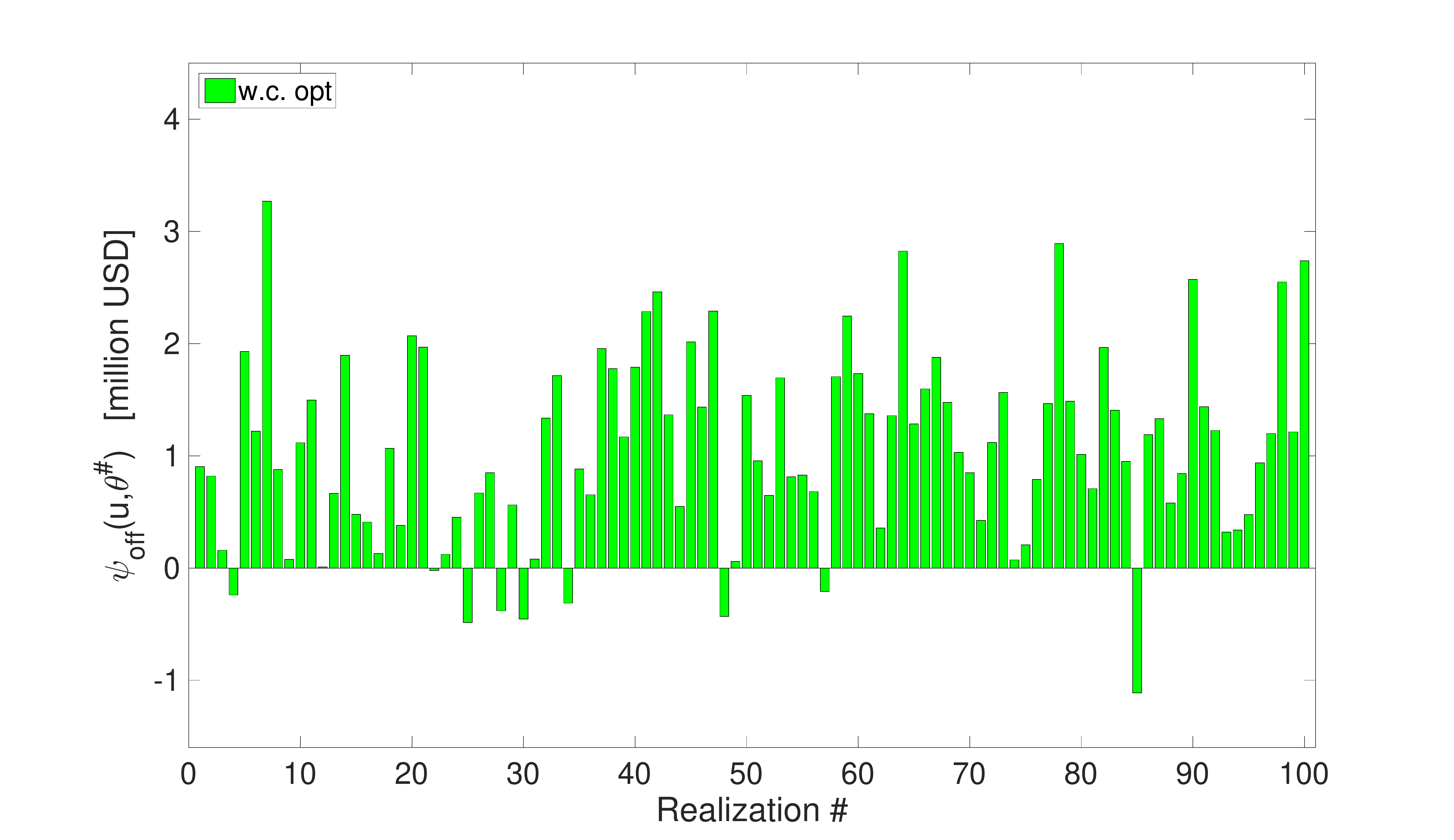}	
		} 
		\subfloat[Profit offset distribution using the worst case control strategy.]{
			\includegraphics[width=0.495\textwidth]{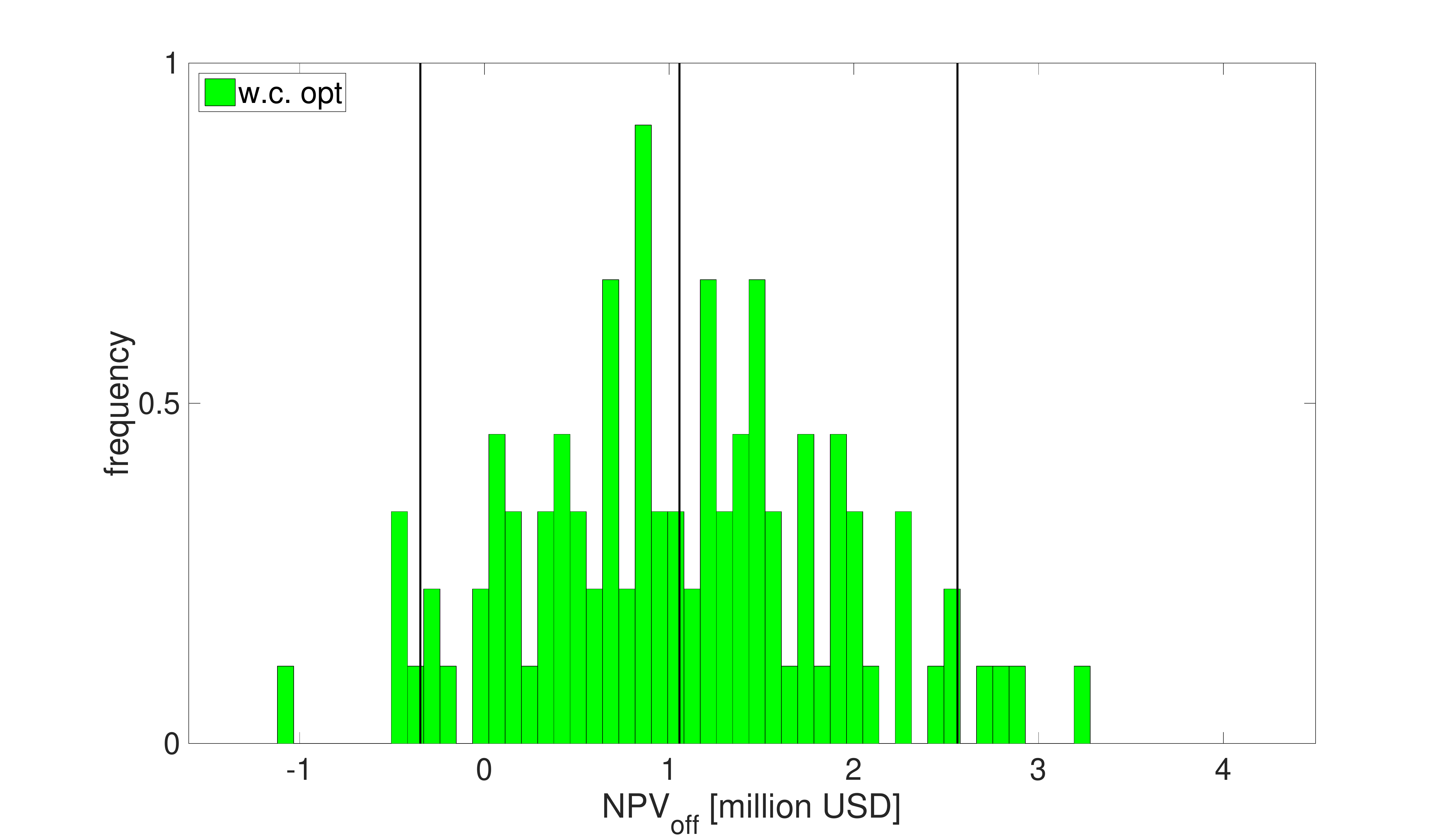}	
		} 
		\\
		\subfloat[Profit offset realizations using the c.s. 30\% control strategy.]{
			\includegraphics[width=0.495\textwidth]{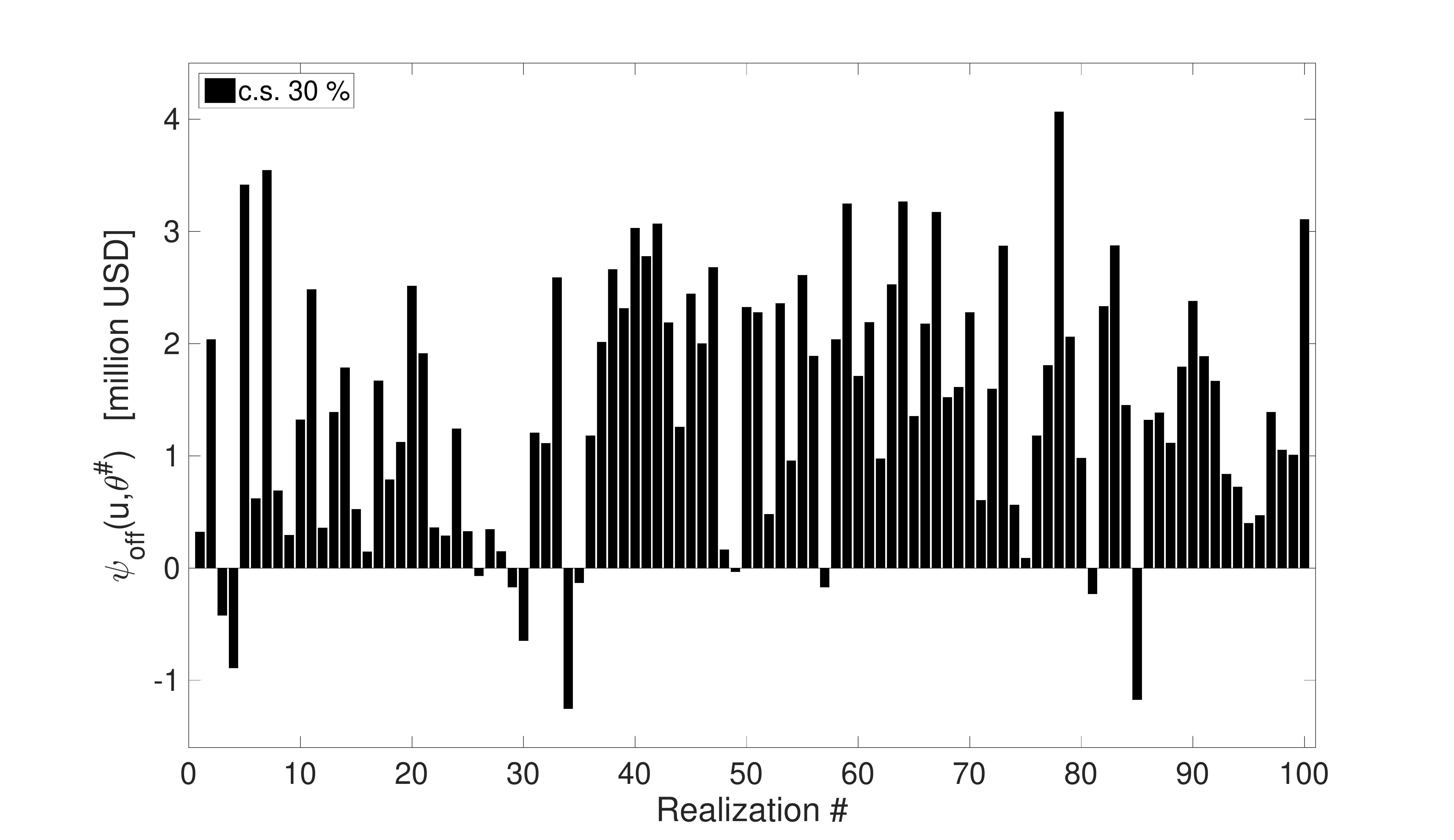}	
		} 
		\subfloat[Profit offset distribution using the c.s. 30\% control strategy.]{
			\includegraphics[width=0.495\textwidth]{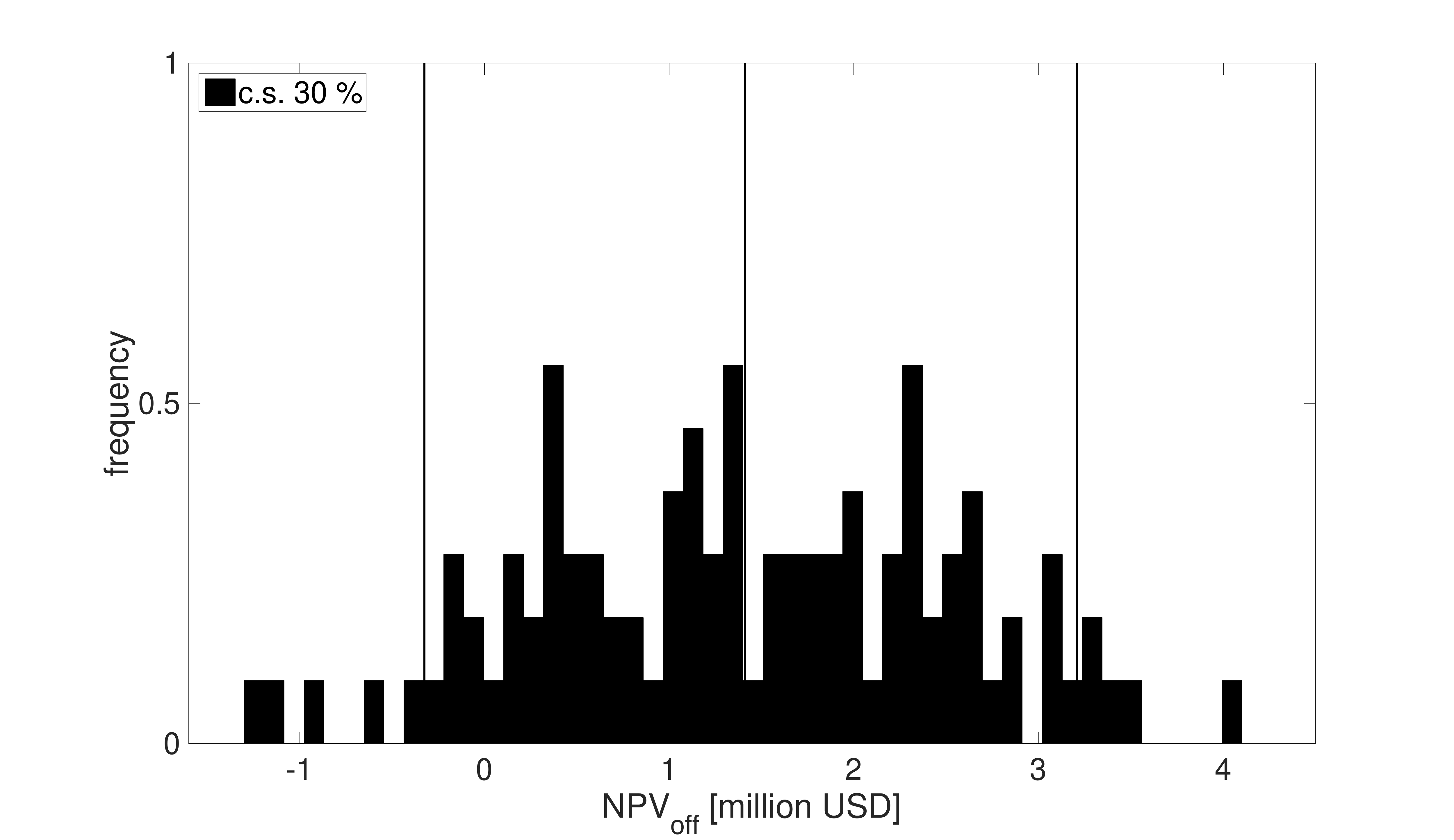}	
		} 
		\\
	 	\subfloat[Profit offset realizations using the RO control strategy.]{
			\includegraphics[width=0.495\textwidth]{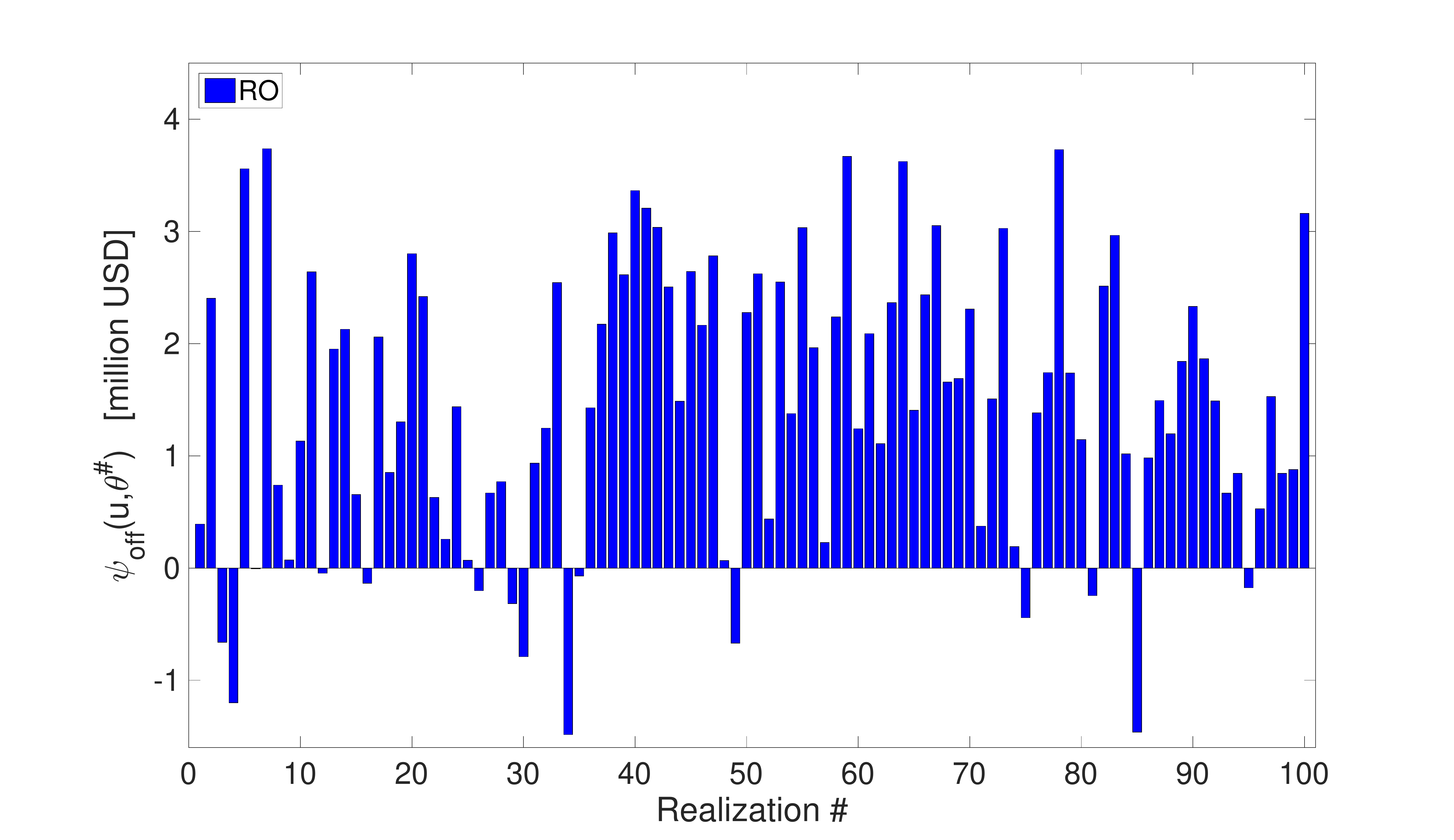}	
		}
		\subfloat[Profit offset distribution using the RO control strategy.]{
			 \includegraphics[width=0.495\textwidth]{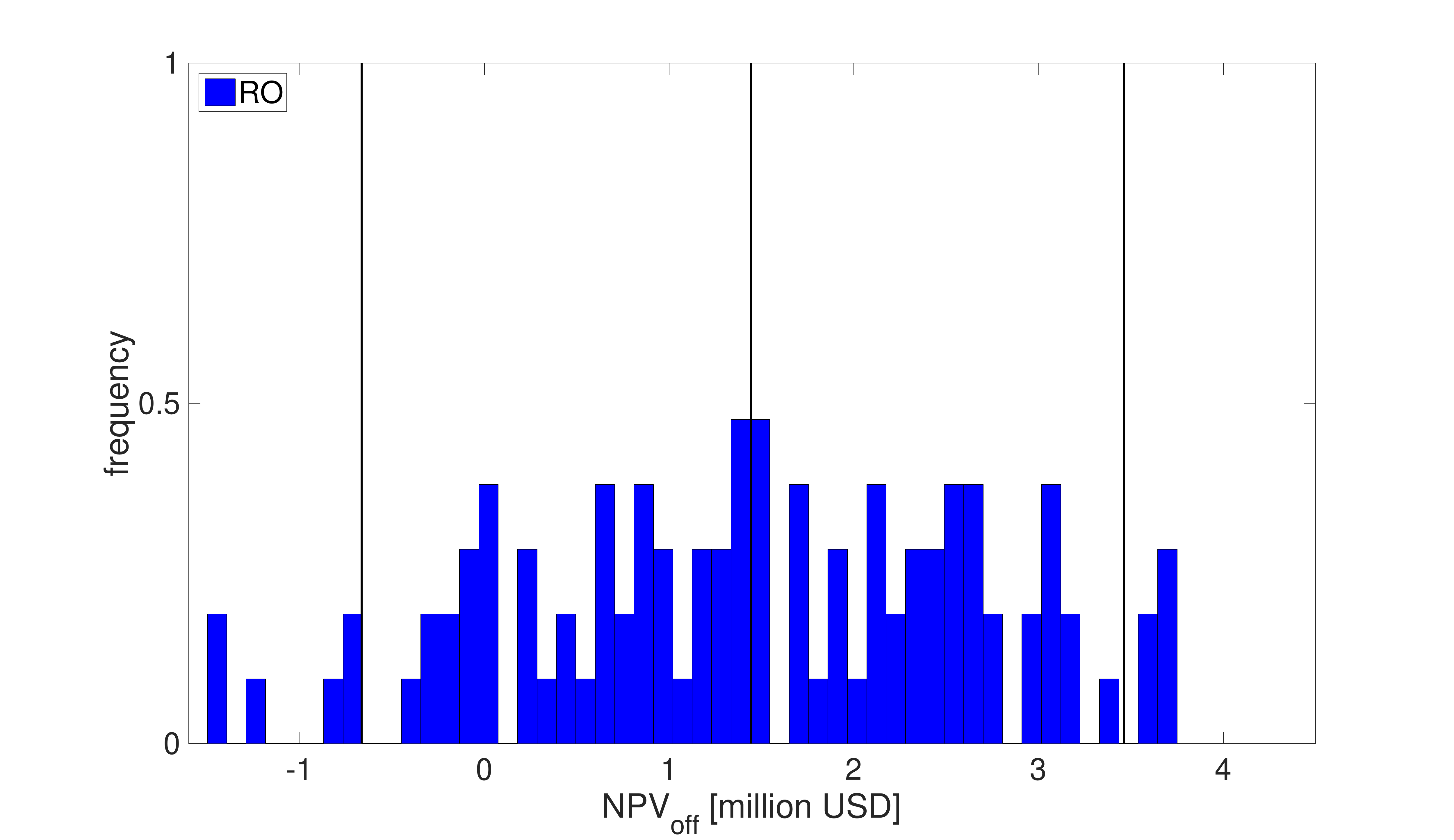}	
		}								
	\end{center}
	\caption{Realizations of the profit offset, $\psi_{off}$, and the corresponding distributions for three different optimal control strategies.}
	\label{fig:OffsetProfitRealizationsAndDistribution}
\end{figure*}

To investigate the risk of the optimizing strategies doing worse than the reactive strategy, we consider the NPV offset computed by (\ref{prodOpt:profExplOff}) for each of the realizations of the optimized case and the reference case. The left panel of Fig. 
\ref{fig:OffsetProfitRealizationsAndDistribution} shows the offset realizations for selected optimized strategies. Evidently, all of the optimized strategies contain realizations that have negative offset NPV, i.e. that do worse than the reactive control strategy. The right panel of Fig. \ref{fig:OffsetProfitRealizationsAndDistribution} illustrates the corresponding profit offset distributions. It is directly visible that the risk that the optimized strategies yield lower profits than the reactive strategy is significant and in the range 8-15\% (see also Table \ref{tb:offSetPsi}).

To handle this situation of profits in relation to a reference strategy, we introduce the profit offset by (\ref{prodOpt:profExplOff}) and an associated risk measure that can be minimized. This allows us to quantify the risk of obtaining profit realizations that are lower than those of the reference strategy, i.e. the reactive strategy. In particular, we say that an optimized strategy provides acceptable risk provided that $\mathcal{R}(\psi_{off}) \leq 0$. For CVaR risk measures this implies that the mean of the $\alpha$-percent worst offset realizations is negative, i.e. that the mean of these optimized strategies is better than the reference strategy.

\begin{figure}[!tb]
	\begin{center}
			\includegraphics[width=0.5\textwidth]{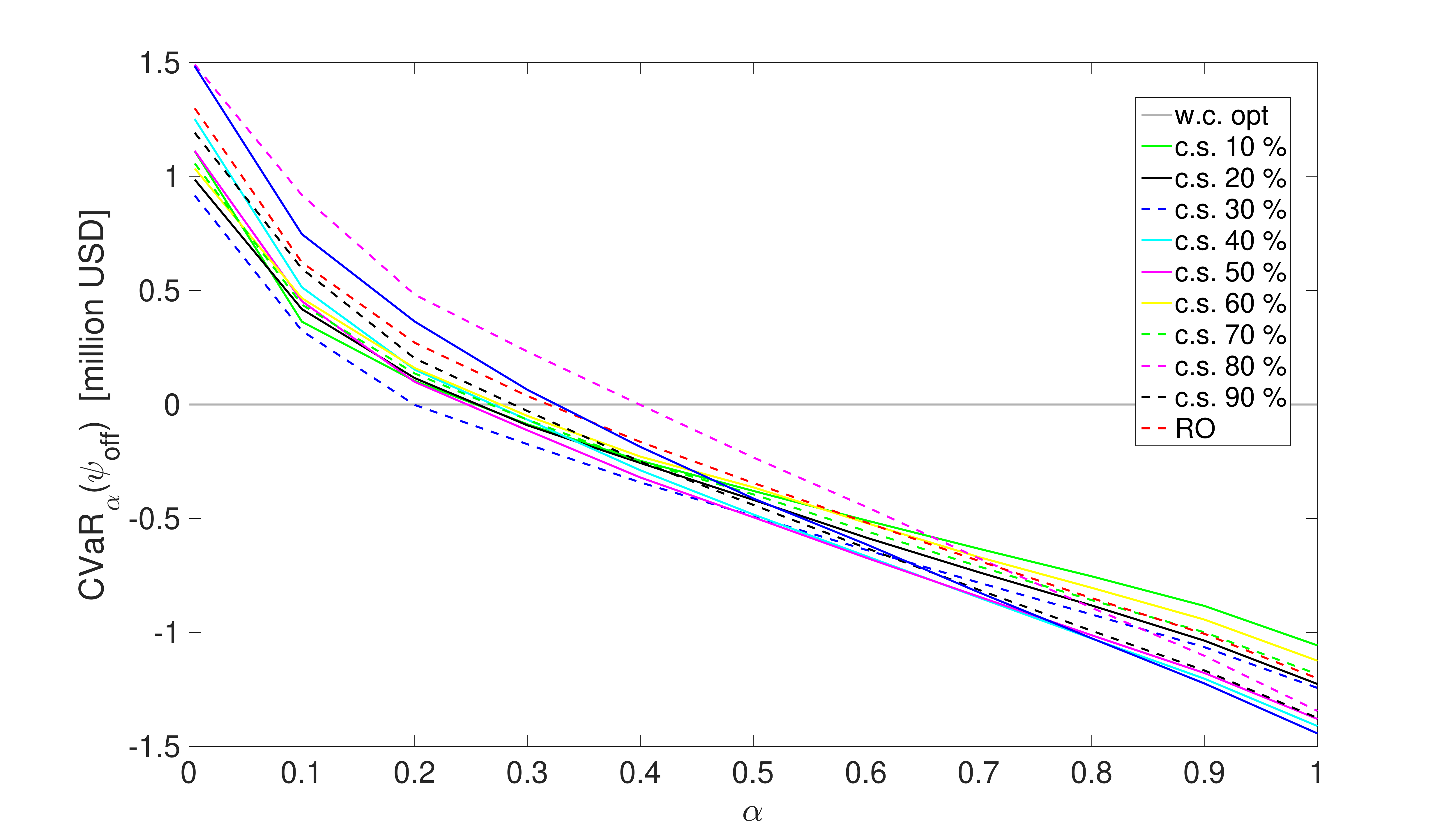}
	\end{center}
	\caption{$\text{CVaR}_{\alpha}(\psi_{off})$ as a function of the risk level, $\alpha$, for different control strategies based on profit risk optimization. Positive values of $\text{CVaR}_{\alpha}$ corresponds to a risk  of doing worse than reactive strategy. For a risk level $\alpha < 20$ in relation to the profit offset, all optimized strategies on average do worse than the reactive strategy, and therefore have a significant (non-negligible) risk of yielding a production profit that is lower than what would be achieved by the reactive strategy.}
	\label{fig:CVarOffSol}
\end{figure}

Fig. \ref{fig:CVarOffSol} shows 11 different CVaR risk measures associated with the NPV offset distribution, $\psi_{off}$. As expected, when the risk level is $\alpha < 20 \%$, none of the optimized control strategies have a $\text{CVaR}_\alpha \leq 0$ for the NPV offset, i.e., all optimized strategies have a positive probability of yielding profits for some realizations that are lower than the reference strategy. 



%

\subsection{Minimization of the profit offset risk}

It is evident that even though the overall risks of the open-loop optimized strategies are much smaller than the reactive strategy, there is still a non-negligible risk that an open-loop optimized strategy for certain realizations may perform worse than the reactive strategy. Therefore, if the reservoir asset manager is more concerned with not doing worse than the reactive strategy, he would select the input profile, $u$, such that the worst profit realization would be maximized, i.e.
\begin{equation}
	\label{eq:MinProfitOffsetRisk1}
		\max_{u \in \mathcal{U}} \inf_{\theta}\left( \psi_{off}(u,\theta) \right) \quad = \quad \max_{u \in \mathcal{U}} \inf_{i=1,\ldots,n_d}\left( \psi_{off}(u,\theta_i) \right).
\end{equation}
The optimization problem \eqref{eq:MinProfitOffsetRisk1} is non-smooth. However, the numerical solution of \eqref{eq:MinProfitOffsetRisk1} is equivalent to the solution of the smooth constrained optimization problem 
\begin{subequations}
	\label{eq:MinProfitOffsetRisk2}
\begin{alignat}{5}
	& \max_{s \in \Real, \, u \in \mathcal{U}} \quad && s, \\
	& s.t. && s \leq \psi_{off}(u,\theta_i), \quad && i=1,\ldots, n_d,
\end{alignat}
\end{subequations}
which may be converted to a minimization problem
\begin{subequations}
		\label{eq:MinProfitOffsetRisk3}
\begin{alignat}{5}
	& \min_{t \in \Real, \, u \in \mathcal{U}} \quad && t, \\
	& s.t. && \psi_{off}(u,\theta_i) + t \geq 0, \quad i=1, \ldots, n_d,
\end{alignat}
\end{subequations}
with an equivalent solution. Furthermore, the solution of \eqref{eq:MinProfitOffsetRisk1} is equivalent to the solution of 
\begin{equation}
	\min_{u \in \mathcal{U}} \left[ -\inf_{i=1,\ldots,n_d} \left( \psi_{off}(u,\theta_i) \right) \right]
	 \, = \, 	\min_{u \in \mathcal{U}} \left[ \text{CVaR}_{\alpha}\left( \psi_{off}(u,\theta) \right) \right], 
\end{equation}
for $\alpha \in [0,p)$. Consequently, we may regard maximization of the worst case offset profit as an offset profit CVaR minimization problem. A natural extension of this interpretation of the worst case offset profit maximization problem is to consider it as a constrained CVaR minimization problem
\begin{subequations}
			\label{eq:MinProfitOffsetRisk4}
\begin{alignat}{5}
& \min_{u \in \mathcal{U}} \quad && \text{CVaR}_{\alpha}(\psi_{off}(u,\theta)), \\
& s.t. && \psi_{off}(u,\theta_i) \geq s,  \quad && i=1,\ldots,n_d, \label{eq:MinProfitOffsetRisk4:constraint}
\end{alignat}
\end{subequations}
for $\alpha \in [0,1]$ and with the parameter $s \in \Real$ denoting the worst acceptable profit offset. To have a feasible solution, $s \leq s^*=-t^*$, where $s^*$ is the solution of \eqref{eq:MinProfitOffsetRisk2} and $t^*$ is the solution of \eqref{eq:MinProfitOffsetRisk3}. The solution, $u = u(\alpha,s)$, is a function of the risk level, $\alpha$, and the worst acceptable offset profit, $s$. 

A number of variations to \eqref{eq:MinProfitOffsetRisk4} exist. One variation is to replace the worst tolerable offset profit \eqref{eq:MinProfitOffsetRisk4:constraint} with a probabilistic constraint 
\begin{subequations}
	\label{eq:MinProfitOffsetRisk5}
	\begin{alignat}{5}
	& \min_{u \in \mathcal{U}} \quad && \text{CVaR}_{\alpha}(\psi_{off}(u,\theta)), \\
	& s.t. && \text{Prob}_\theta[ \psi_{off}(u,\theta) < 0 ] \leq \beta,  \label{eq:MinProfitOffsetRisk5:constraint}
	\end{alignat}
\end{subequations}
for $\alpha \in [0,1]$ and $\beta \in [0,1]$. The constraint \eqref{eq:MinProfitOffsetRisk5:constraint} expresses that the probability of having a negative offset profit, i.e. doing worse than the reference strategy, should be less than $\beta$. Obviously, the feasibility of such a constraint depends on the value of $\beta$ and the specific reservoir being studied. 

\begin{figure*}[!tb]
 	\begin{center}
 		\subfloat[Strip chart of the profit offset realizations.]{
 			\includegraphics[width=0.5\textwidth]{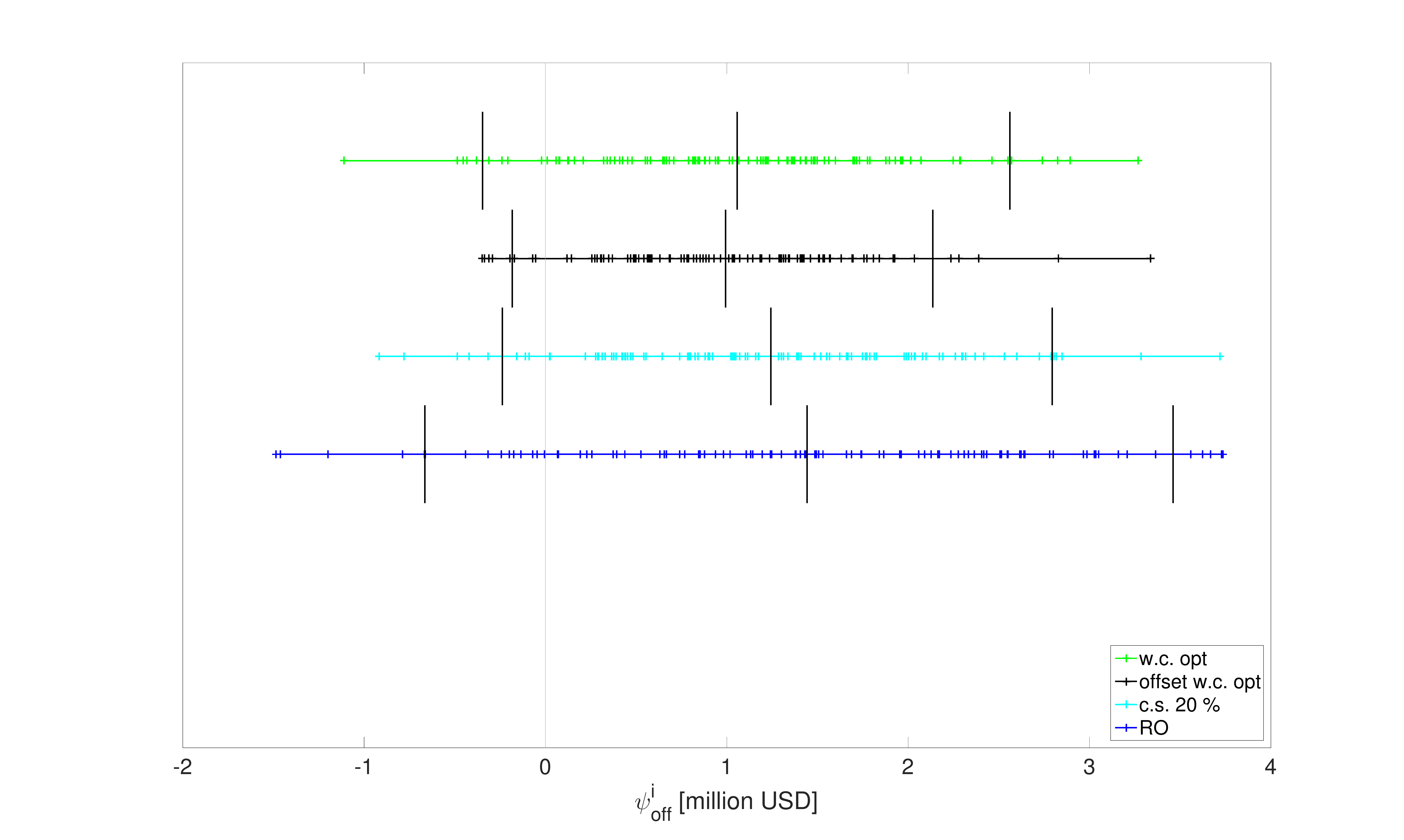}
 			\label{fig:CVarOffDistrWor}}
 		\subfloat[$\text{CVaR}_{\alpha}(\psi_{off})$ as a function of the risk level, $\alpha$.]{
 			\includegraphics[width=0.5\textwidth]{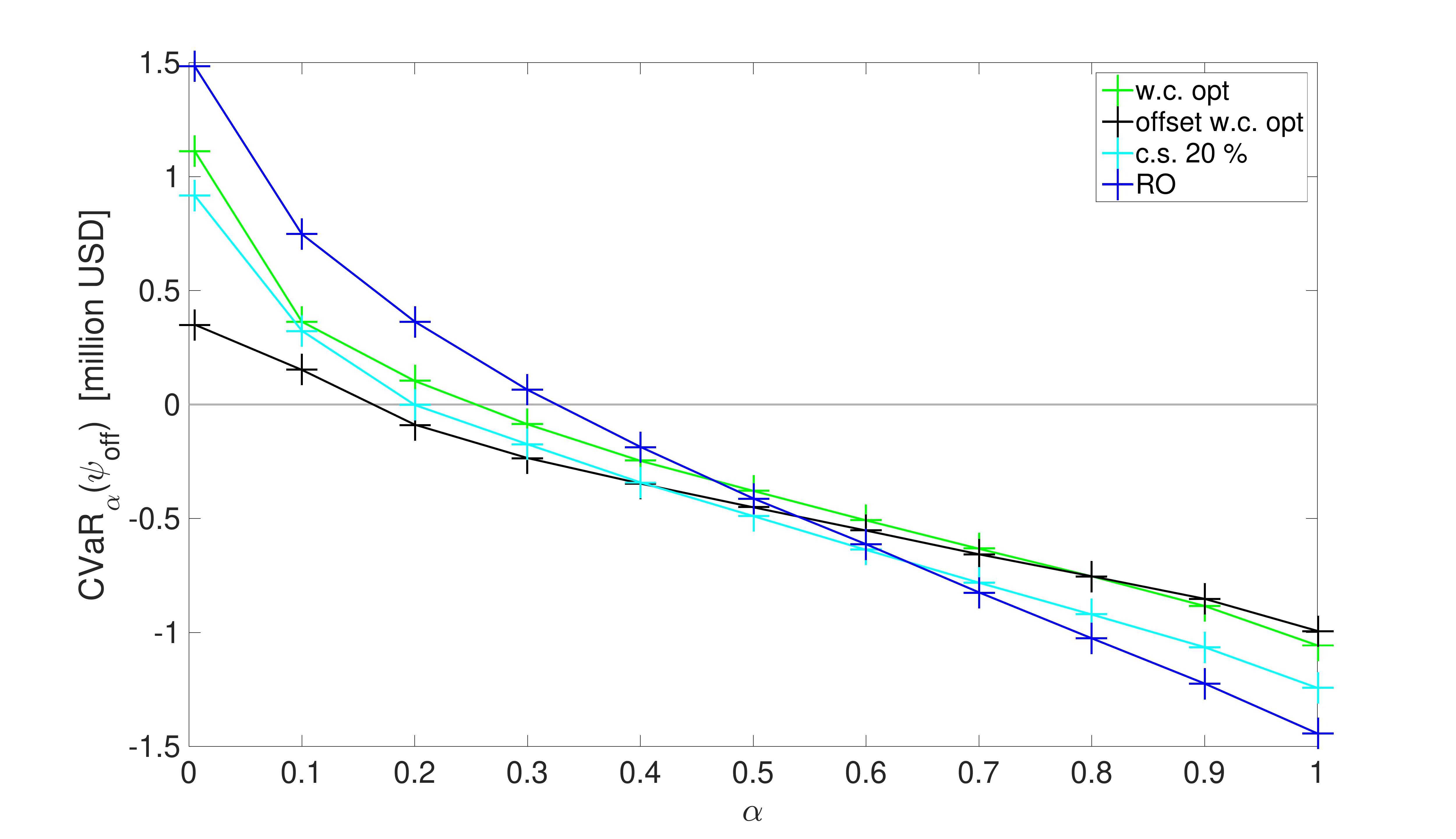}  
 			\label{fig:CVarOffSol1}
 		}
 	\end{center}
 	\caption{The profit offset realization and the CVaR$_\alpha$ for selected control strategies. (a) Strip charts of the NPV offset distribution for selected optimized control strategies. The black vertical lines indicate the 5th percentile, the mean, and the 95th percentile of the profit offset distribution. The offset worst case optimization strategy has less risk of doing worse than the reactive strategy but also lower expected profit compared to the strategies optimizing profit. (b) Plot of $\text{CVaR}_{\alpha}(\psi_{off})$ as a function of the risk level, $\alpha$, for different optimized control strategies. For low risk levels, $\alpha < 0.2$, the optimized control strategies all have a risk of doing worse than the reactive strategy. The offset worst case optimization strategy has the lowest risk of doing worse than the reactive strategy.   } 
 	\label{fig:offEnsNPVsAll}
 \end{figure*}

In this paper, we solve \eqref{eq:MinProfitOffsetRisk4}  for $\alpha \in [0,p)$ and $s=s^*=-t^*$, i.e. we solve \eqref{eq:MinProfitOffsetRisk1}. This control strategy is called offset worst case optimization and is denoted "offset w.c. opt." (see Table \ref{tb:namingConvContr}). Fig. \ref{fig:offEnsNPVsAll} illustrates the profit offset realizations and the CVaR$_{\alpha}$ at different risk levels, $\alpha$, for the offset worst case optimization strategy and selected profit optimization strategies, i.e. the worst case optimization strategy (w.c. opt), the CVaR$_{20\%}$ optimization strategy (c.s. 20\%), and the robust optimization strategy (RO). Fig. \ref{fig:CVarOffDistrWor} shows the NPV offset distributions and demonstrates the offset worst case optimization strategy produces some realization that do worse than the reactive strategy. However, the realizations indicate that the loss and the number of occurrences with negative offset are smaller for the worst case profit offset optimization compared to the strategies based on profit optimization. This is confirmed by the CVaR plot of the profit offset in Fig. \ref{fig:CVarOffSol1}. The reduced risk of the offset worst case optimization strategy in relation to the reactive straegy comes as the price of less expected profit than the strategies minimizing risk related to the profit and not the profit offset.

 Table \ref{tb:offSetPsi} provides key performance indicators for the control strategies based on profit optimization and the profit offset worst case  optimization strategy. The mean profit offset column, $\text{E}_{\theta}(\psi_{off})$, shows that all optimized strategies do better then the reactive strategy on average. The column with the worst profit offset realization, $\inf(\psi_{off})$, shows that all optimized strategies risk producing realizations with negative profit offset, i.e. risk doing worse than the reactive strategy for some realizations. The profit offset worst case optimization strategy has a better worst case profit offset realization but a lower expected profit realization than the strategies minimizing risk related to profit. The fourth to the sixth column of Table \ref{tb:offSetPsi} report the probability of having a negative profit offset, $\beta = \text{Prob}[\psi_{off} < 0]$, the average profit offset of the negative profit offsets, $\text{E}_{\theta}[\psi_{off}|\psi_{off}<0]$, and the average profit offset of the positive profit offsets, $\text{E}_{\theta}[\psi_{off}|\psi_{off}\geq0]$. This implies that the control strategy produced by the offset worst case optimization (off w.c. opt.) has 8\% chance of a negative profit offset that has an average offset profit of -0.22 mio USD, and a 92\% chance of a positive profit offset with an average value of 1.10 mio USD. The offset worst case optimization yields the lowest probability of negative profit offset and these negative profit offsets have the highest average value among the control strategies investigated. However, the offset worst case optimization also has the lowest average offset profit among the positive offset profits. This implies that offset worst case optimization produces a control strategy with the lowest risk, in a certain sense, of doing worse than the reactive strategy, but it also has the least ability to improve the reactive strategy. The last block of columns in Table \ref{tb:offSetPsi} reports the key performance indicators related to the profit. From an overall profit perspective, other optimization based strategies exist (i.e. w.c opt., c.s. 10\%, c.s. 20\%) that have both a higher worst case profit, a higher average profit, and a lower risk (higher -CVaR$_{30\%}$) than the strategy produced by the profit offset worst case optimization. In summary, Table \ref{tb:offSetPsi} demonstrates that we can produce a strategy, based on profit offset worst case optimization, that provides the least risk of doing worse than the reactive strategy, have a high probability (92\%) of doing better, and that in term of worst case profit, expected profit, and risk (CVaR$_{30\%}$) is within 1-2.5\% of the best of these measures for the optimized strategies. 
 
 Inspection of Table \ref{tb:offSetPsi} provides the impression that all optimized strategies provide similar performance. However, the optimized strategies perform somewhat better than the reactive strategy, even though the reactive strategy uses feedback, while the optimized strategies are based on open-loop optimization, i.e. no feedback. This impression is also confirmed by Fig. \ref{fig:BarPlot}. It shows that the largest improvement comes between the optimized strategies and the reactive strategy and not between the different optimized strategies.

\begin{table*}[!tb]
	\begin{center}
		\caption{Key performance indicators for the NPV offset distribution.}\label{tb:offSetPsi}
		\scalebox{0.7}{
			\begin{tabular}{l|cc|ccc|ccccccr}
				\hline\\[-8pt]
				Control 
				
				& $\text{E}_{\theta}(\psi_{off})$  
				& $\inf(\psi_{off})$  
				& $\beta := \text{Prob}[\psi_{off} < 0]$  
				& $\text{E}_{\theta}[\psi_{off} | \psi_{off}  < 0]$ 
				& $\text{E}_{\theta}[\psi_{off} \big| \psi_{off} \geq 0]$   
				& $\inf (\psi)$    & $\text{E}_{\theta}(\psi)$ & $\sigma_{\theta}(\psi)$ & $5\%$ perc. & $95\%$ perc.  & -CVaR$_{30\%}(\psi)$  
				\\
				strategy & 	$10^6$ USD	 				& $10^6$ USD	 & & $10^6$ USD   & $10^6$ USD 	& $10^6$ USD 	& $10^6$ USD 	& $10^6$ USD  	& $10^6$ USD 	& $10^6$ USD 	& $10^6$ USD \\\hline
				w.c. opt.   		  &	1.06 & -1.11 &  9\% 	&  -0.41 & 1.20     	 &  42.94             &	45.43 &  1.43    & 43.34 & 47.83 & 43.70      \\
				c.s. $10\%$     			 & 1.23 & -0.99  &  10\%    &  -0.42 &  1.41     &  42.14        &  45.60 &  1.44 & 43.47 & 48.13 & 43.88       \\
				c.s. $20\%$   			 & 1.24 & -0.92  &  8\%  &  -0.41 &  1.39      &  41.94        &  45.62 &  1.44   & 43.47 & 48.12 & 43.89      \\
				c.s. $30\%$   			 & 1.41 & -1.25      &  11\%  &  -0.47 &  1.64      	 &  41.65        &  45.79 &  1.48  & 43.15 & 48.34 & 44.03   \\
				c.s. $40\%$   			 & 1.38 & -1.11    &  8\%      &  -0.58 &  1.55       &  41.70       &  45.76 &  1.49 & 43.13 & 48.30 & 43.98          \\
				c.s. $50\%$   			& 1.12 & -1.04    &  11\%        &  -0.43 & 1.32    	 &  41.48       &  45.50 & 1.51 & 43.13 & 47.98 & 43.67           \\
				c.s. $60\%$   			 & 1.18 & -1.06   &  10\%        &  -0.44 &  1.36     	 &  41.49        &  45.56 &  1.53 & 43.10 & 48.11 & 43.72     \\
				c.s. $70\%$   			 & 1.34 & -1.49     &  15\%      &  -0.66 &  1.70      	 &  40.94       &  45.72 &  1.64 & 42.79 & 48.31 & 43.76         \\
				c.s. $80\%$   			  & 1.37 & -1.19     &  9\%     &  -0.66 &  1.58    	 &  41.46      &  45.75 &  1.55  & 43.01 & 48.26 & 43.90           \\
				c.s. $90\%$   			  & 1.20 & -1.30    &  12\%      &  -0.55 &  1.44      	 &  41.23        &  45.58 &  1.54 & 42.87 & 48.12 & 43.73    \\
				RO		   			 & 1.44 & -1.48     &  15\%      &  -0.53 &  1.79      &  41.45    &  45.82 &  1.58 & 43.06 & 48.45 & 43.95          \\ \hline
				offset w.c. opt.    &	0.99 & -0.35 	 &  8\% 	 & -0.22 & 1.10    & 41.94 & 45.37 & 1.55 & 42.64 & 47.73 & 43.49 \\  \hline
				ref & 0.00 & 0.00 & 0\% & 0.00 & 0.00 & 40.60 & 44.38 & 1.57 & 41.59 & 46.57 & 42.46 \\ \hline
			\end{tabular} }
		\end{center}
	\end{table*}

\section{Conclusions}\label{conclusions}
This paper explores the concept of risk minimization in life-cycle oil production optimization. In this context, we propose to use the axioms of coherence and aversion as a systematic approach to characterize proper risk measures. As a specific example of a proper measure, we consider conditional value-at-risk, $\text{CVaR}_{\alpha}$, at different risk levels, $\alpha$. By a numerical case study, we investigate the ability of $\text{CVaR}_{\alpha}$ to minimize the risk of profit losses.  As a benchmark reference strategy, representing real-world best practices, we use reactive control. By minimizing $\text{CVaR}_{\alpha}$ over an ensemble of 100 permeability realizations, we show that for every risk level, $\alpha$, there exists optimized  strategies that yield a lower risk than the risk obtained using a reactive control strategy. However, the results also show, that in spite of the overall lower risk, the optimized strategies are associated with a significant risk of yielding low profit outcomes relative to reactive control. To mitigate this risk of very low profit realizations, we introduce a method that seeks to optimize the worst-case NPV offset value in relation to a reference strategy, i.e. the reactive strategy. Minimizing the risk \textit{relative} to a reference strategy is novel compared to existing methods available in the open oil literature. 
Using the offset CVaR approach, we significantly reduce the risk of low profit outcomes in relation to the current best practice, i.e. the reactive strategy. In particular, relative to the widely used RO strategy, we find an optimized strategy that manages to halve both the probability of having a low profit outcome and the actually profit loss of such a low outcome. As a minor drawback, we do not find an optimized strategy with zero probability of yielding lower profit realizations than the reactive strategy. This is most likely because the reference strategy incorporates valuable feedback via reactive control. The optimized strategies studied in this paper are all so-called open-loop strategies that do not employ feedback. Future work will seek to overcome this issue by combining the optimization procedure with feedback e.g. using a receding horizon implementation of combined data assimilation and optimization.

\section*{Acknowledgement}
This project is financially supported by The Danish Advanced Technology Foundation (OPTION; 63-2013-3).

\appendix
\section{Nomenclature} \label{Sec:app}
\begin{table}[h!]
\centering
\scalebox{0.66}{
\begin{tabular}{ll} 
\hline
 Symbol & Description \\
\hline
 $c_{ref}$   & deterministic reference payoff \\
 		 $\text{CVaR}_{\alpha}$  &  conditional value at risk at confidence level $\alpha$\\ 
 		  $\text{E}_{\theta} (\cdot)$  &   expected value      \\ 
 		 KKT  &   Karush-Kuhn-Tucker conditions   \\ 
 		 $n_d$  &   number of possible scenarios     \\ 
 		 $n_u$  &   dimension of the control vector     \\ 
 		 $\text{Prob}[\cdot]$     &  probability of an event \\ 
 		 RO     &   Robust Optimization strategy     \\ 
 		 $u$  &   control vector     \\ 
 		 $u_{opt}$  &  optimized control strategy     \\ 
 		 $u_{ref}$  &  reference control strategy     \\
		\hline 
		 $\mathcal{U}$  &   set of feasible controls     \\ 
		  $p$  &   probability of a scenario  under equiprobability assumptions  \\ 
		  $p_i$  &   probability of i-th scenario\\ 
		  $\mathcal{R}$  &   risk measure \\
 		 $\mathcal{R}_{total}$  &   total risk measure\\ 
 		 $\text{VaR}_{\alpha}$  &   value at risk at confidence level $\alpha$\\
		\hline
$\beta$  &   probability of getting a profit  lower than the reference strategy\\ 
 		$\theta$  &  random permeability field vector     \\ 
 		$\theta_i$  &  i-th scenario for permeability field \\ 
 		$\Theta$  &   uncertainty space of the permeability field      \\ 
 		$\Theta_d$  &   discretized uncertainty space of the permeability field      \\ 
 		$\sigma^2 (\cdot)$  &   variance       \\ 
 		$\psi$    & random profit distribution  \\ 
 		$\psi^i$    &   profit outcome of the i-th scenario   \\ 
 		$\psi_{off}$    &  random profit offset distribution  \\
 		$\psi^i_{off}$    &   profit offset outcome of the i-th scenario   \\  \hline
\end{tabular} \label{fig:paraetwell} 
} 
\end{table}

\bibliography{ReferencesACAP}

\end{document}